\documentclass[12pt]{article}
\usepackage{amsthm, amsmath, amssymb, enumerate}
\usepackage{fullpage}
\usepackage[colorlinks=true]{hyperref}

\begin{document}
\title{Explicit Kodaira--Spencer map over Shimura Curves}
\author{Xinyi Yuan}
\maketitle

\theoremstyle{plain}
\newtheorem{thm}{Theorem}[section]
\newtheorem{theorem}[thm]{Theorem}
\newtheorem{cor}[thm]{Corollary}
\newtheorem{corollary}[thm]{Corollary}
\newtheorem{lem}[thm]{Lemma}
\newtheorem{lemma}[thm]{Lemma}
\newtheorem{pro}[thm]{Proposition}
\newtheorem{proposition}[thm]{Proposition}
\newtheorem{prop}[thm]{Proposition}
\newtheorem{definition}[thm]{Definition}
\newtheorem{assumption}[thm]{Assumption}

\theoremstyle{remark} 
\newtheorem{remark}[thm]{Remark}
\newtheorem{example}[thm]{Example}
\newtheorem{remarks}[thm]{Remarks}
\newtheorem{problem}[thm]{Problem}
\newtheorem{exercise}[thm]{Exercise}
\newtheorem{situation}[thm]{Situation}

\numberwithin{equation}{subsection}

\newcommand{\ZZ}{\mathbb{Z}}
\newcommand{\CC}{\mathbb{C}}
\newcommand{\QQ}{\mathbb{Q}}
\newcommand{\RR}{\mathbb{R}}
\newcommand{\HH}{\mathcal{H}}     

\newcommand{\ad}{\mathrm{ad}}            
\newcommand{\NT}{\mathrm{NT}}         
\newcommand{\nonsplit}{\mathrm{nonsplit}}         
\newcommand{\Pet}{\mathrm{Pet}}         
\newcommand{\Fal}{\mathrm{Fal}}         

\newcommand{\cs}{{\mathrm{cs}}}         

\newcommand{\ZZn}{\mathbb{Z}[1/n]}         
\newcommand{\ZZN}{\mathbb{Z}[1/N]}         

\newcommand{\pair}[1]{\langle {#1} \rangle}
\newcommand{\wpair}[1]{\left\{{#1}\right\}}
\newcommand{\wh}{\widehat}
\newcommand{\wt}{\widetilde}

\newcommand\Spf{\mathrm{Spf}}

\newcommand{\lra}{{\longrightarrow}}

\newcommand{\matrixx}[4]
{\left( \begin{array}{cc}
  #1 &  #2  \\
  #3 &  #4  \\
 \end{array}\right)}        


\newcommand{\BA}{{\mathbb {A}}}
\newcommand{\BB}{{\mathbb {B}}}
\newcommand{\BC}{{\mathbb {C}}}
\newcommand{\BD}{{\mathbb {D}}}
\newcommand{\BE}{{\mathbb {E}}}
\newcommand{\BF}{{\mathbb {F}}}
\newcommand{\BG}{{\mathbb {G}}}
\newcommand{\BH}{{\mathbb {H}}}
\newcommand{\BI}{{\mathbb {I}}}
\newcommand{\BJ}{{\mathbb {J}}}
\newcommand{\BK}{{\mathbb {K}}}
\newcommand{\BL}{{\mathbb {L}}}
\newcommand{\BM}{{\mathbb {M}}}
\newcommand{\BN}{{\mathbb {N}}}
\newcommand{\BO}{{\mathbb {O}}}
\newcommand{\BP}{{\mathbb {P}}}
\newcommand{\BQ}{{\mathbb {Q}}}
\newcommand{\BR}{{\mathbb {R}}}
\newcommand{\BS}{{\mathbb {S}}}
\newcommand{\BT}{{\mathbb {T}}}
\newcommand{\BU}{{\mathbb {U}}}
\newcommand{\BV}{{\mathbb {V}}}
\newcommand{\BW}{{\mathbb {W}}}
\newcommand{\BX}{{\mathbb {X}}}
\newcommand{\BY}{{\mathbb {Y}}}
\newcommand{\BZ}{{\mathbb {Z}}}

\newcommand{\CA}{{\mathcal {A}}}
\newcommand{\CB}{{\mathcal {B}}}
\newcommand{\CD}{{\mathcal{D}}}
\newcommand{\CE}{{\mathcal {E}}}
\newcommand{\CF}{{\mathcal {F}}}
\newcommand{\CG}{{\mathcal {G}}}
\newcommand{\CH}{{\mathcal {H}}}
\newcommand{\CI}{{\mathcal {I}}}
\newcommand{\CJ}{{\mathcal {J}}}
\newcommand{\CK}{{\mathcal {K}}}
\newcommand{\CL}{{\mathcal {L}}}
\newcommand{\CM}{{\mathcal {M}}}
\newcommand{\CN}{{\mathcal {N}}}
\newcommand{\CO}{{\mathcal {O}}}
\newcommand{\CP}{{\mathcal {P}}}
\newcommand{\CQ}{{\mathcal {Q}}}
\newcommand{\CR }{{\mathcal {R}}}
\newcommand{\CS}{{\mathcal {S}}}
\newcommand{\CT}{{\mathcal {T}}}
\newcommand{\CU}{{\mathcal {U}}}
\newcommand{\CV}{{\mathcal {V}}}
\newcommand{\CW}{{\mathcal {W}}}
\newcommand{\CX}{{\mathcal {X}}}
\newcommand{\CY}{{\mathcal {Y}}}
\newcommand{\CZ}{{\mathcal {Z}}}

\newcommand{\ab}{{\mathrm{ab}}}
\newcommand{\Ad}{{\mathrm{Ad}}}
\newcommand{\an}{{\mathrm{an}}}
\newcommand{\Aut}{{\mathrm{Aut}}}

\newcommand{\Br}{{\mathrm{Br}}}
\newcommand{\bs}{\backslash}
\newcommand{\bbs}{\|\cdot\|}

\newcommand{\Ch}{{\mathrm{Ch}}}
\newcommand{\cod}{{\mathrm{cod}}}
\newcommand{\cont}{{\mathrm{cont}}}
\newcommand{\cl}{{\mathrm{cl}}}
\newcommand{\criso}{{\mathrm{criso}}}

\newcommand{\dR}{{\mathrm{dR}}}
\newcommand{\disc}{{\mathrm{disc}}}
\newcommand{\Div}{{\mathrm{Div}}}
\renewcommand{\div}{{\mathrm{div}}}

\newcommand{\Eis}{{\mathrm{Eis}}}
\newcommand{\End}{{\mathrm{End}}}

\newcommand{\Frob}{{\mathrm{Frob}}}

\newcommand{\Gal}{{\mathrm{Gal}}}
\newcommand{\GL}{{\mathrm{GL}}}
\newcommand{\GO}{{\mathrm{GO}}}
\newcommand{\GSO}{{\mathrm{GSO}}}
\newcommand{\GSp}{{\mathrm{GSp}}}
\newcommand{\GSpin}{{\mathrm{GSpin}}}
\newcommand{\GU}{{\mathrm{GU}}}
\newcommand{\BGU}{{\mathbb{GU}}}

\newcommand{\Hom}{{\mathrm{Hom}}}
\newcommand{\Hol}{{\mathrm{Hol}}}
\newcommand{\HC}{{\mathrm{HC}}}

\renewcommand{\Im}{{\mathrm{Im}}}
\newcommand{\Ind}{{\mathrm{Ind}}}
\newcommand{\inv}{{\mathrm{inv}}}
\newcommand{\Isom}{{\mathrm{Isom}}}

\newcommand{\Jac}{{\mathrm{Jac}}}
\newcommand{\JL}{{\mathrm{JL}}}

\newcommand{\Ker}{{\mathrm{Ker}}}
\newcommand{\KS}{{\mathrm{KS}}}

\newcommand{\Lie}{{\mathrm{Lie}}}

\newcommand{\new}{{\mathrm{new}}}
\newcommand{\NS}{{\mathrm{NS}}}

\newcommand{\ord}{{\mathrm{ord}}}
\newcommand{\ol}{\overline}

\newcommand{\rank}{{\mathrm{rank}}}

\newcommand{\PGL}{{\mathrm{PGL}}}
\newcommand{\PSL}{{\mathrm{PSL}}}
\newcommand{\Pic}{\mathrm{Pic}}
\newcommand{\Prep}{\mathrm{Prep}}
\newcommand{\Proj}{\mathrm{Proj}}

\newcommand{\Picc}{\mathcal{P}ic}

\renewcommand{\Re}{{\mathrm{Re}}}
\newcommand{\red}{{\mathrm{red}}}
\newcommand{\sm}{{\mathrm{sm}}}
\newcommand{\sing}{{\mathrm{sing}}}
\newcommand{\reg}{{\mathrm{reg}}}

\newcommand{\tor}{{\mathrm{tor}}}
\newcommand{\tr}{{\mathrm{tr}}}

\newcommand{\ur}{{\mathrm{ur}}}

\newcommand{\vol}{{\mathrm{vol}}}

\newcommand{\ds}{\displaystyle}

\tableofcontents

\section{Introduction}

The goal of this paper is to explicitly compute the Kodaira--Spencer map for a quaternionic Shimura curve over $\QQ$ and its effect on the metrics of the Hodge bundle. 
It is easy to define the relevant maps canonically, but it is not easy to get the explicit constants due to various abstract identifications. In fact, there are many wrong explicit results in the literature.

Let $B$ be an indefinite quaternion algebra over $\QQ$. 
Denote by $\tr:B\to \QQ$ and $q:B\to \QQ$ the reduced trace and the reduced norm.
Denote the main involution by $\iota:B\to B, \ \beta\mapsto \beta^\iota$.
Denote by $d_B$ the discriminant of $B$, which is the product of the prime numbers 
$p$ at which $B$ is a division algebra.

Fix a maximal order $O_B$ of $B$ throughout this paper.
The main involution stabilize $O_B$ by the formula $\beta^\iota=\tr(\beta)-\beta$.
Let $U=\prod_p U_p$ be an open subgroup of $\wh O_B^\times$. 
Let $n=n(U)$ be the product of primes $p$ such that $U_p$ is not maximal, i.e. $U_p\neq O_{B,p}^\times$. 

Associated to these data, there is a Shimura curve $X_U$ over $\CC$, and the Shimura curve has a canonical integral model $\CX_U$ over $\ZZ[1/n]$.
In fact, $\CX_U$ is a stack over $\ZZ[1/n]$ such that for any $\ZZ[1/n]$-scheme $S$, $\CX_U(S)$ is the category of triples 
$(A, i,\bar\eta)$ as follows:
\begin{enumerate}[(1)]
\item $A$ is an abelian scheme of relative dimension 2 over $S$;
\item $i:O_B\to \End_S(A)$ is a ring homomorphism satisfying the determinant condition  
$$
\det(\beta|_{\Lie(A/S)})=q(\beta),\quad\forall \beta\in O_B;
$$
\item $\bar\eta$ is a $U$-level structure on $A$. 
\end{enumerate}
Then $\CX_U$ is a flat regular Deligne--Mumford stack over $\ZZ[1/n]$, but it is actually a regular scheme if $U$ is sufficiently small.
We refer to \S\ref{sec Shimura} for more details.

The moduli structure does not indicate any polarization, but it turns out that a polarization is automatic.
Fix an element $\mu\in O_B$ such that $\mu^2=-d_B$ throughout this paper. This gives a positive involution 
on $B$ by 
$\beta^*=\mu^{-1} \beta^\iota \mu$.
Then for any triple $(A, i,\bar\eta)$ over $S$ as above, there exists a unique principle polarization
$\lambda:A\to A^t$ whose Rosati involution on $\End_S(A)$ is compatible with the positive involution $*$ on $O_B$ via $i$.

Denote by $\pi:\CA\to \CX_U$ the universal abelian scheme, and by $\lambda:\CA\to \CA^t$ the universal principle polarization (depending on $\mu$). 
Denote by $\epsilon:\CX_U\to \CA$ and $\epsilon^t:\CX_U\to \CA^t$ the identity sections.

We have the relative differential sheaves $\Omega_{\CA/\ZZ[1/n]}$, $\Omega_{\CA/\CX_U}$ and $\Omega_{\CX_U/\ZZ[1/n]}$, the relative dualizing sheaves $\omega_{\CA/\CX_U}$ and $\omega_{\CX_U/\ZZ[1/n]}$, and the relative tangent sheaves $T_{\CA/\CX_U}$ and $T_{\CX_U/\ZZ[1/n]}$ (of derivations). 
We further have the Lie algebra 
$$
\Lie(\CA):=\epsilon^* T_{\CA/\CX_U}\simeq \pi_*T_{\CA/\CX_U},
$$
and the \emph{Hodge bundles}
$$
\underline\Omega_{\CA}:=\epsilon^* \Omega_{\CA/\CX_U}\simeq \pi_*\Omega_{\CA/\CX_U},\qquad
\underline\omega_{\CA}:=\epsilon^* \omega_{\CA/\CX_U}\simeq \pi_*\omega_{\CA/\CX_U}.
$$
We have easy canonical isomorphisms
$$
\underline\Omega_{\CA}\simeq \Lie(\CA)^\vee, \qquad
\underline\omega_{\CA}\simeq\det\underline\Omega_{\CA}.
$$
These definitions work for Deligne--Mumford stacks via \'etale descent.

Endow the Hodge bundle $\underline\omega_{\CA}$ over $\CX_U$ with the 
\emph{Faltings metric}
$\|\cdot\|_\Fal$ as follows.
For any point $x\in \CX_U(\CC)$, and any section $\alpha\in \underline\omega_{\CA}(x)\simeq \Gamma(\CA_x,\omega_{\CA_x/\CC})$, the Faltings metric is defined by
$$
\|\alpha\|_\Fal^2:=\frac{1}{(2\pi )^2}\left|\int_{\CA_x(\CC)} \alpha\wedge\overline\alpha \right|.
$$
Here $\CA_x$ is the fiber of $\CA$ above $x$, which is a complex abelian surface, 
$\alpha$ is viewed as a holomorphic 2-form over $\CA_x$ via the canonical isomorphism
$\underline\omega_{\CA}(x)\simeq \Gamma(\CA_x,\omega_{\CA_x/\CC})$.

Note that $\CX_U$ is regular, so the dualizing sheaf $\omega_{\CX_U/\ZZ[1/n]}$ is a line bundle over $\CX_U$.
Endow $\omega_{\CX_U/\ZZ[1/n]}$ over $\CX_U$ with the 
\emph{Petersson metric}
$\|\cdot\|_\Pet$ as follows.
Via the complex uniformization of $\CH$ to every connected component of $\CX_U(\CC)$, 
it suffices to define a metric on $\omega_{\CH/\CC}$ that descends to $\omega_{\CX_U(\CC)/\CC}$. 
Denote by $\tau$ the usual coordinate function of $\CH\subset\CC$. 
Then we have a canonical identification $\omega_{\CH/\CC}=\CO_{\CH} d\tau$.
The Petersson metric is defined by 
$$
\|d\tau\|_\Pet:=2\,\Im(\tau). 
$$

Finally, we are ready to introduce the Kodaira--Spencer map. 
Start with the exact sequence
$$
0\lra \pi^*\Omega_{\CX_U/\ZZ[1/n]}\lra \Omega_{\CA/\ZZ[1/n]}\lra \Omega_{\CA/\CX_U}\lra 0.
$$
Apply derived functors of $\pi_*$.
It gives a connecting morphism 
$$
\phi_0:\pi_*\Omega_{\CA/\CX_U} \lra R^1\pi_*(\pi^*\Omega_{\CX_U/\ZZ[1/n]}).
$$
This is the \emph{Kodaira--Spencer map}. 
There are canonical isomorphisms 
$$R^1\pi_*(\pi^*\Omega_{\CX_U/\ZZ[1/n]}) \lra 
R^1\pi_*\CO_{\CA}\otimes \Omega_{\CX_U/\ZZ[1/n]}\lra 
\Lie(\CA^t)\otimes \Omega_{\CX_U/\ZZ[1/n]}
\lra \underline\Omega_{\CA^t}^\vee\otimes \Omega_{\CX_U/\ZZ[1/n]},$$ 
where $\CA^t\to\CX_U$ denotes the dual abelian scheme of $\CA\to \CX_U$.
Then the {Kodaira--Spencer map} is also written as 
$$
\phi_1:\underline\Omega_{\CA} \lra \underline\Omega_{\CA^t}^\vee\otimes \Omega_{\CX_U/\ZZ[1/n]}.
$$

There is a canonical isomorphism $\Omega_{\CX_U/\ZZ[1/n]}\to \omega_{\CX_U/\ZZ[1/n]}$ over the smooth locus $\CX_U^\sm$ of $\CX_U$ over $\ZZ[1/n]$, and this isomorphism extends to a morphism over the $\CX_U$ due to the fact that $\CX_U$ is regular and $\CX_U\setminus\CX_U^\sm$ has codimension 2.
Then the Kodaira--Spencer map induces a morphism 
$$
\phi_2:\underline\Omega_{\CA} \lra \underline\Omega_{\CA^t}^\vee\otimes \omega_{\CX_U/\ZZ[1/n]}.
$$
Taking determinants, we get a morphism
$$
\psi_1:\underline\omega_{\CA} \lra \underline\omega_{\CA^t}^\vee\otimes \omega_{\CX_U/\ZZ[1/n]}^{\otimes 2}.
$$
This further induces a canonical morphism 
$$
\psi_2:\underline\omega_{\CA} \otimes \underline\omega_{\CA^t}\lra \omega_{\CX_U/\ZZ[1/n]}^{\otimes 2}.
$$
Finally, via the principal polarization $\lambda:\CA\to \CA^t$, it induces a canonical morphism 
$$
\psi_3:\underline\omega_{\CA}^{\otimes 2}\lra \omega_{\CX_U/\ZZ[1/n]}^{\otimes 2}.
$$ 

The following is the main result of this paper.
\begin{thm}\label{KS}
The canonical map 
$\psi_3:\underline\omega_{\CA}^{\otimes 2}\to \omega_{\CX_U/\ZZ[1/n]}^{\otimes 2}$ is injective, and its image is the subsheaf $d_B\,\omega_{\CX_U/\ZZ[1/n]}^{\otimes 2}$ of $\omega_{\CX_U/\ZZ[1/n]}^{\otimes 2}$. 
Moreover, under $\psi_3$, we have $\|\cdot\|_\Fal^2=\|\cdot\|_\Pet^2$.
\end{thm}

From the result, we see that if $B$ is a division algebra, then there is no morphism $\underline\omega_{\CA}\to \omega_{\CX_U/\ZZ[1/n]}$
whose square is equal to $\psi_3:\underline\omega_{\CA}^{\otimes 2}\to \omega_{\CX_U/\ZZ[1/n]}^{\otimes 2}$.

Note that our results agree with \cite[Theorem 3.7, Theorem 4.10]{YZ}, but do not agree with \cite[Proposition 3.2, Lemma 3.3]{KRY1}.

If $B$ is a division algebra and
$U=\wh O_B^\times$ is maximal, then $\CX_U$ is proper over $\ZZ$, and we can compute arithmetic intersection numbers of the above hermitian line bundles. 
As a consequence of Theorem \ref{KS}, 
$$
\frac{\widehat\deg(\hat c_1(\omega_{\CX_U/\ZZ}, \|\cdot\|_\Pet)^2)}{2\deg(\omega_{\CX_{U,\QQ}/\QQ})}
=\frac{\widehat\deg(\hat c_1(\underline\omega_{\CA}, \|\cdot\|_\Fal)^2)}{2\deg(\underline\omega_{\CA,\QQ})}
+\frac12 \log d_B.
$$
This is precisely stated in Theorem \ref{height}, which also includes arithmetic intersection numbers over the coarse moduli schemes.
A similar result for modular curves (in the case $B=M_2(\QQ)$) is in Theorem \ref{height2}.
Then the main theorem of Yuan \cite{Yua}, a formula on the modular height of quaternionic Shimura curves over totally real fields, 
is compatible with the formulas of K\"uhn \cite{Kuh} and Kudla--Rapoport--Yang \cite{KRY2} over $\QQ$.
This the main motivation of this paper.

\begin{remark}
In this paper, we only consider integral models $\CX_U$ over $\ZZ[1/n]$, but one can extend them to proper integral models $\CX_U^*$ over $\ZZ$ using Drinfel'd level structures. Then Theorem \ref{KS} for $\CX_{\wh O_B^\times}$ implies a similar result for $\CX_U^*$ by considering the morphism $\CX_U^*\to \CX_{\wh O_B^\times}$.
\end{remark}

This paper is organized as follows. 
In \S\ref{sec Shimura}, we introduce more on Shimura curves and their coarse moduli schemes.
In  \S\ref{sec image}, we prove the first statement of Theorem \ref{KS}, which follows from deformation theory with a lengthy extra argument at places $p|d_B$. 
In \S\ref{sec metric}, we prove the second statement of Theorem \ref{KS}, which follows from a lengthy explicit computation of the Kodaira--Spencer map in terms of complex uniformization.
In \S\ref{sec modular curve}, we extend the main results of this paper to 
the Deligne--Rapoport compactification of modular curves.

\subsubsection*{Acknowledgment}
The author would like to thank Kai-Wen Lan and Congling Qiu for many helpful communications.

\section{Shimura curves and modular curves} \label{sec Shimura}

In this section, we review some basics on integral models of Shimura curves. 
While our main result in Theorem \ref{KS} is described in term of moduli stacks, we will also consider its consequence in terms of coarse moduli schemes here.

\subsection{Shimura curves}

Let us review integral models of Shimura curves following the expositions of \cite{Mil, BC, Buz}. 
We will repeat some terminology of the previous section, but give more details in the following.

Let $B$ be an indefinite quaternion algebra over $\QQ$. 
Denote by $\tr:B\to \QQ$ and $q:B\to \QQ$ the reduced trace and the reduced norm.
Denote the main involution by $\iota:B\to B, \ \beta\mapsto \beta^\iota$.
Denote by $d_B$ the discriminant of $B$, which is the product of prime numbers 
$p$ at which $B$ is a division algebra.

Fix a maximal order $O_B$ of $B$ throughout this paper.
The main involution stabilize $O_B$ by the formula $\beta^\iota=\tr(\beta)-\beta$.
Denote $\wh O_{B}=O_{B}\otimes_\ZZ \wh\ZZ$ and $O_{B,p}=O_{B}\otimes_\ZZ \ZZ_p$ for any prime $p$.
Let $U=\prod_p U_p$ be an open subgroup of $\wh O_B^\times$. 
Let 
$$n=n(U)=\prod_{p:\,U_p\neq O_{B,p}^\times} p$$ 
be the product of primes $p$ such that $U_p$ is not maximal. 

By assumption, $U$ contains $U(N)$ for some positive integer $N$, 
where 
$$U(N)=\prod_{p\nmid N}  O_{B,p}^\times \times \prod_{p| N} (1+N O_{B,p})^\times$$
is the principal open subgroup of $\wh O_{B}^\times$.

Fix an isomorphism $\sigma:B\otimes_\ZZ\RR\simeq M_2(\RR)$ throughout this paper.
Then we have a complex Shimura curve
$$X_U=B^\times \backslash \CH^{\pm}\times B_{\BA_f}^\times/U
=B_+^\times \backslash \CH\times B_{\BA_f}^\times/U.$$
Here $B^\times$ acts on $\CH^\pm$ via $\sigma$, and 
$$B_+^\times=\{\gamma\in B^\times: q(\gamma)>0\}.$$
Note that $X_U$ is understood to be an orbifold, but it is actually a smooth complex curve if the quotient process is free, which happens if $U$ contains $U(N)$ for some $N\geq3$.

The Shimura curve $X_U$ has a canonical integral model $\CX_U$ over $\ZZ[1/n]$ for $n=n(U)$.
In fact, $\CX_U$ is a stack over $\ZZ[1/n]$ such that for any $\ZZ[1/n]$-scheme $S$, $\CX_U(S)$ is the category of triples 
$(A, i,\bar\eta)$ as follows:
\begin{enumerate}[(1)]
\item $A$ is an abelian scheme of pure relative dimension 2 over $S$;
\item $i:O_B\to \End_S(A)$ is a ring homomorphism satisfying the determinant condition  
$$
\det(\beta|_{\Lie(A/S)})=q(\beta),\quad\forall \beta\in O_B;
$$
\item $\bar\eta$ is a $U$-level structure on $A$. 
\end{enumerate}
The $U$-level structure is defined as follows. Note that $U$ contains $U(N)$ for some positive integer $N$ dividing a power of $n$. 
Consider pairs $(T,\eta)$, where $T\to S$ is an \'etale cover of $S$, and 
$\eta: O_B/N O_B\to A(T)[N]$ is an $O_B$-linear map  satisfying $\eta_2=\eta_1\circ u$ for some $u\in U$ acting on $O_B/N O_B$, where  $\eta_j: O_B/N O_B\to A[N]$  is the composition of $\eta$ with the map $p_j^*:A(T)[N]\to A(T\times_ST)[N]$ induced by the projection $p_j: T\times_ST\to T$ for $j=1,2$. 
Two pairs $(T,\eta)$ and $(T',\eta')$ are said to be $U$-equivalent if there is an \'etale cover $V\to S$ refining both $T\to S$ and $T'\to S$ such that $\eta'_{V}=\eta_{V}\circ u$ for some $u\in U$ acting on $O_B/N O_B$, where
$\eta_{V}: O_B/N O_B\to A(V)[N]$ and $\eta'_{V}: O_B/N O_B\to A(V)[N]$ are the maps induced by $\eta$ and $\eta'$ respectively.
Finally, a $U$-level structure on $A$ is an equivalence class of pairs $(T,\eta)$.

As mentioned in the introduction, there exists a unique principal polarization compatible with the moduli structure.
In fact, fix an element $\mu\in O_B$ such that $\mu^2=-d_B$ throughout this paper. This gives a positive involution 
$$
*:B\lra B,\quad
\beta^*=\mu^{-1} \beta^\iota \mu, \quad \beta\in B.
$$
The involution stabilizes $O_B$ as a consequence of the same property at every finite place.
Then for any triple $(A, i,\bar\eta)$ over $S$ as above, there exists a unique principle polarization
$\lambda:A\to A^t$ whose Rosati involution on $\End_S(A)$ is compatible with the positive involution $*$ on $O_B$ via $i$.
This result is a consequence of  \cite[Lemma 1.1]{Mil}, \cite[III, Proposition 1.5, Proposition 3.3]{BC}, and \cite[\S8,\S11]{Bou}. 

It is well-known that $\CX_U$ is a regular Deligne--Mumford stack, flat and semistable over $\ZZ[1/n]$, and smooth outside $d_B$. 
Moreover, it is actually a scheme if $U$ contains $U(N)$ for some $N\geq3$.
If $B$ is a division algebra, which is our major concern, then
$\CX_U$ is further proper over $\ZZ[1/n]$. 

As a Deligne--Mumford stack, $\CX_U$ has an \'etale cover by schemes, so 
most terminologies and properties of schemes can be transferred to $\CX_U$ via \'etale descent. 
Explicitly in the current situation, for any open subgroup $U'\subset U$, the morphism $\CX_{U'}\to \CX_{U}$ is \'etale.
Let $U'$ be sufficiently small so that $\CX_{U'}$ is a scheme. 
Then all such $\CX_{U'}$ form an \'etale cover of $\CX_U$ by schemes. 
In fact, if $U$ is maximal, we can always find two such $\CX_{U'}$ to cover $\CX_U$; if $U$ is not maximal, we can find a single $\CX_{U'}$ to cover $\CX_U$.
Because of that, in our treatment, we can usually reduce the problem for the stack $\CX_U$ to that for the scheme $\CX_{U'}$.

\subsection{Coarse moduli schemes} \label{sec coarse}

Now we review the notion of coarse moduli schemes and Hodge bundles over them. 
Our exposition is a special case of \cite[\S4]{YZ}.

Denote by $\CX_U^\cs$ the \emph{coarse moduli scheme} of $\CX_U$ over $\ZZn$. 
This can be constructed explicitly as follows.
Let $U'\subset U$ be a normal and open subgroup such that $\CX_{U'}$ is a (regular) scheme over $\ZZ[1/n(U')]$.
The finite group $U/U'$ acts on the regular scheme $\CX_{U'}$.
Denote by 
$\CX_{U'}/U=\CX_{U'}/(U/U')$ the quotient \emph{scheme} over $\ZZ[1/n(U')]$.
Recall that the quotient scheme is locally defined by the spectrum of the ring of invariant functions.
If $U''$ is another subgroup of that type, then $\CX_{U'}/U$ and $\CX_{U''}/U$ are canonically isomorphic over  $\ZZ[1/(n(U')n(U''))]$. 
Thus we can glue all these $\CX_{U'}$ to form a scheme over $\ZZn$, and the result is the coarse moduli scheme $\CX_U^\cs$.  

By definition, $\CX_U^\cs$ is a 2-dimensional normal scheme, flat and quasi-projective over $\ZZn$. 
The generic fiber $\CX^\cs_{U,\QQ}$ of $\CX^\cs_{U}$ is a smooth curve over $\QQ$ by the quotient process.
If $B$ is a division algebra, then $\CX_U^\cs$ is projective over $\ZZn$. 
Crucial to intersection theory, the quotient process further implies that $\CX_U^\cs$ is 
\emph{$\QQ$-factorial}, i.e. any Weil divisor over $\CX_U^\cs$ has a positive multiple which is a Cartier divisor.

To prove the last statement, it suffices to prove it for an effective Weil divisor $\CD$ over $\CX_{U'}/U$ for $i=1,2$. 
The pullback of $\CD$ to $\CX_{U'}$ is locally defined by a regular function $f$ on $\CX_{U'}$.
The key is that $\mathrm{N}(f)=\prod_{\gamma\in U/U'} \gamma(f)$ is a local regular function on $\CX_{U'}/U$, and the divisor  $\div(\mathrm{N}(f))=[U:U'] \CD$. 
This proves the {$\QQ$-factorial} property.

To introduce the Hodge bundle over the coarse moduli scheme, we first introduce the category of $\QQ$-line bundles in an abstract setting. 
For any Deligne--Mumford stack $S$, denote by $\Picc(S)$ the \emph{category} of line bundles over $S$, in which the objects are line bundles (or equivalently invertible sheaves) over $S$, and the morphisms are morphisms of coherent sheaves.  
Denote by $\Picc(S)_\QQ$ the \emph{category} of $\QQ$-line bundles over $S$, in which the objects are pairs $(a,\CL)$ (or just written as $\CL^{\otimes a}$)
with $a\in \QQ$ and $\CL\in\Picc(S)$, and the morphism of two such objects is defined to be
$$\Hom(\CL^{\otimes a},\CL'^{\otimes a'}):=\varinjlim_m
 \Hom(\CL^{\otimes (am)}, \CL'^{\otimes (a'm)}),$$
where $m$ runs through positive integers such that $am$ and $a'm$ are both integers, 
and ``$\Hom$'' on the right-hand side represents isomorphisms of integral line bundles. 
For the direct system, for any $m|k$, there is a transition map 
$$\Hom(\CL^{\otimes (am)}, \CL'^{\otimes (a'm)})\lra \Hom(\CL^{\otimes (ak)}, \CL'^{\otimes (a'k)})$$ 
locally given by taking $(k/m)$-th power. 
A \emph{section} of a $\QQ$-bundle $\CL^{\otimes a}$ is an element of $\Hom(\CO_S, \CL^{\otimes a})$.
We can also introduce rational sections and their divisors similarly. 

Denote by $(\CX_U^\cs)^\reg$ the regular locus of $\CX_U^\cs$.
Since $\CX_U^\cs$ is normal, the complement of $(\CX_U^\cs)^\reg$ in $\CX_U^\cs$ is 0-dimensional. 
Then the relative dualizing sheaf $\omega_{(\CX_U^\cs)^\reg/\ZZn}$
extends to a $\QQ$-line bundle over $\CX_U$, which can be checked by the passage between line bundles and divisors. 
Denote this extension by $\omega_{\CX_U^\cs/\ZZn}$, and call it the 
\emph{relative dualizing sheaf} of $\CX_U^\cs$ over $\ZZn$.
It is unique up to unique isomorphisms.

Finally, the \emph{Hodge bundle} $\CL_U$ is a $\QQ$-line bundle over $\CX_U^\cs$ defined by
$$
\CL_U:= \omega_{\CX_{U}^\cs/\ZZn} \otimes 
\left(\bigotimes_{P\in \CX^\cs_{U,\QQ}} \CO_{\CX_U^\cs}(\CP)^{\otimes (1-e_P^{-1})}\right).
$$
Here $\omega_{\CX_{U}/\ZZn}$ is the relative dualizing sheaf defined above, 
the summation is through closed points $P$ on the generic fiber $\CX^\cs_{U,\QQ}$ of $\CX^\cs_{U}$ over $\QQ$, 
$\CP$ is the Zariski closure of $P$ in $\CX_U^\cs$, and $e_P$ is the ramification index  of $P$ in the map $\CX_{U'}\to \CX_{U}^\cs$ 
for any normal and open subgroup $U'\subset U$ such that $\CX_{U'}$ is a scheme.
One can check that $e_P$ does not depend on the choice of $U'$, and that 
$e_P$ is also equal to the ramification index by the uniformization map from $\CH$ to connected components of $\CX_U^\cs(\CC)$.

If $\CX_{U}$ is already a scheme,
then we simply have $\CX_{U}^\cs=\CX_{U}$; if furthermore 
$B$ is a division algebra, then 
$\CL_{U}= \omega_{\CX_{U}/\ZZn}.$
The following result justifies the definition involving the ramification indices.

\begin{lem} \label{compatibility}
\begin{itemize}
\item[(1)]
Let $U'\subset U$ be a normal and open subgroup such that $\CX_{U'}$ is a scheme.
Let $\pi_{U',U}:\CX_{U'}\to \CX_{U}^\cs$ be the natural morphism.
Then there is a canonical isomorphism
$$
\pi_{U',U}^*\CL_{U} \lra \CL_{U'}
$$
of $\QQ$-line bundles over $\CX_{U'}$.
\item[(2)]
Denote by $f:\CX_{U}\to \CX_{U}^\cs$ the canonical morphism.
Then there is a canonical isomorphism
$$
f^*\CL_{U} \lra \omega_{\CX_U/\ZZ}
$$
of $\QQ$-line bundles over $\CX_{U}$.
\end{itemize}
\end{lem}

\begin{proof}

Note that (1) implies (2) by descent via the \'etale morphism $\CX_{U'}\to \CX_U$. 

For (1), we first have a canonical isomorphism
$\pi_{U',U}^*\CL_{U,\QQ} \to \CL_{U',\QQ}$
of the generic fibers in the setting of the classical Hurwitz formula. 
This gives a section of 
$\pi_{U',U}^*\CL_{U,\QQ}^\vee \otimes \CL_{U',\QQ}$, and thus a rational section $s$ of 
$\pi_{U',U}^*\CL_{U}^\vee \otimes \CL_{U'}$ over $\CX_{U'}$.
It suffices to prove that the $\QQ$-divisor $\CD=\div(s)$ is 0 in $\Div(\CX_{U'})_\QQ$. 
We already know that $\CD$ is 0 on the generic fiber. 
It suffices to prove that the support of $\CD$ 
does not contain any irreducible component of a closed fiber of $\CX_{U'}$ over 
$\ZZ[1/n(U')]$. 
Equivalently, it suffices to prove that $\pi_{U',U}:\CX_{U'}\to \CX_{U}^\cs$ is \'etale at 
the generic point $\xi$ of any irreducible component of the fiber of $\CX_{U'}$ above any prime $p$.
Assume that $U'_p=U_p$ is maximal; otherwise, there is nothing to prove.
We will separate the situation into two cases depending on whether $p$ is split in $B$. 

We first treat the case that $p$ is split in $B$. 
Assume that $\xi$ corresponds to a moduli triple $(A, i,\bar\eta)$ over $\xi$. 
To prove result, it suffices to prove that the stabilizer of $U/U'$ on $(A, i,\bar\eta)$ is equal to $(\{\pm 1\}U')/U'$, so that the quotient process is free at $\xi$ (up to $\{\pm 1\}$). 
For that, it suffices to prove $\Aut(A, i)=\{\pm 1\}$, and thus it suffices to prove $\End(A, i)\simeq \ZZ$. 
Denote by $C$ the Zariski closure of $\xi$ in $\CX_{U'}$, which is a curve over $\BF_p$.
Let $K$ be an imaginary quadratic field which can be embedded in $B$. 
This is achieved by requiring $K$ to be nonsplit at all prime factors of $d_B$. 
We further assume that $p$ is split in $K$.  
There is an ordinary elliptic curve $E$ over $\ol\BF_p$ with $\End(E) \simeq O_K$, which can be obtained as the reduction of an elliptic curve over a number field with complex multiplication by $O_K$.
As $K$ can be embedded in $B$, we have $B\otimes_\QQ K\simeq M_2(K)$, so 
$B$ can be embedded in $M_2(K)$.    
As $\End(E^2)\simeq M_2(O_K)$, there is an embedding $i': O_B\to \End(E^2)$. 
With a level structure,  $(E^2,i', \bar\eta')$ defines a point of $\CX_{U', \ol \BF_p}$. 
By adjusting $\bar\eta'$, we can assume that $(E^2,i', \bar\eta')$ actually lies in $C$. 
Then $(E^2,i', \bar\eta')$ is a specialization of $(A, i,\bar\eta)$, which gives an embedding
$\End(A,i)\hookrightarrow \End(E^2,i')$. 
By $B\otimes_\QQ K\simeq M_2(K)$, the centralizer of $B$ in $M_2(K)$ is isomorphic to $K$, so $\End(E^2,i')\otimes_\QQ K$ is isomorphic to $K$. 
It follows that $\End(A,i)$ has an embedding to $K$.
Varying $K$, we see that $\End(A,i)$ has an embedding to $\QQ$, and thus $\End(A,i)\simeq \ZZ$. This proves the split case. 

Now we treat the case that $p$ is nonsplit in $B$. 
We need the $p$-adic uniformization theorem of \v Cerednik and Drinfel'd (cf. \cite{BC}). 
Many ingredients of our treatment is already in \cite[\S4.6]{Yua}. 
Recall that the uniformization gives an isomorphism  
$$ \wh\CX_{U', \ZZ_p^\ur} =D^{\times}\bs 
\wh\Omega_{\ZZ_p^\ur} \times B_{\BA_f}^\times/{U'}.$$
Here we denote
$$\wh\CX_{U', \ZZ_p^\ur}=\wh\CX_{{U'}}\times_{\Spf\, \ZZ_p} \Spf\, \ZZ_p^\ur,\quad
\wh\Omega_{\ZZ_p^\ur}=\wh\Omega \times_{\Spf\, \ZZ_p} \Spf\, \ZZ_p^\ur,$$
where
$\wh\CX_{U'}$ denotes the formal completion of $\CX_{U'}$ along the special fiber above $p$, $\ZZ_p^\ur$ denotes the completion of the maximal unramified extension of $\ZZ_p$, and $\wh \Omega$ is Deligne's formal integral model of Drinfe'ld's upper half plane $\Omega$ over $\ZZ_p$. 
The quaternion algebra $D$ over $\QQ$ is obtained by changing the invariants of $B$ at $p, \infty$. 
The group $D_p^\times\cong \GL
_2(\QQ_p)$ acts on $\wh\Omega$ by some fractional linear transformation, and on 
$D_p^\times/U_p'\cong \ZZ$ via translation by $v_p\circ q=v_p\circ\det$.

The special fiber of 
$\wh\Omega$ is a union of infinitely many copies of $\BP^1_{\BF_p}$ indexed by homethety classes of lattices in $\QQ_p^2$. 
Let $C_0$ be the irreducible component of the special fiber of $\wh \Omega$ corresponding to the standard lattice $\ZZ_p^2$. 
By this interpretation, the stabilizer of $\GL_2(\QQ_p)$ on $C_0$ is 
$\QQ_p^\times \GL_2(\ZZ_p)$.
Moreover, $\QQ_p^\times \GL_2(\ZZ_p)$ acts on $C_0=\BP^1_{\BF_p}$ by the natural action via the quotient $\QQ_p^\times \GL_2(\ZZ_p) \to \PGL_2(\ZZ_p) \to \PGL_2(\BF_p)$. 
Then the common stabilizer of $\GL_2(\QQ_p)$ on all points of $C_0$ is 
$\QQ_p^\times (1+p M_2(\ZZ_p))$.

For the maximal compact subgroup $U$, a further quotient process gives
$$ \wh\CX_{U, \ZZ_p^\ur}^\cs =D^{\times}\bs 
\wh\Omega_{\ZZ_p^\ur} \times B_{\BA_f}^\times/{U}
=(\{\pm1\}\bs D^{\times})\bs 
\wh\Omega_{\ZZ_p^\ur} \times B_{\BA_f}^\times/{U}.$$
To prove that every irreducible component of the special fiber of $\wh\CX_{U', \ZZ_p^\ur}$ is mapped separably to its image in $\wh\CX_{U, \ZZ_p^\ur}^\cs$, it suffices to prove that 
every irreducible component $\wt C$ of the special fiber of 
$\wh\Omega_{\ZZ_p^\ur} \times B_{\BA_f}^\times/{U}$ is mapped separably to its image in $\wh\CX_{U, \ZZ_p^\ur}^\cs$. 
Then it suffices to prove that the common stabilizer of $D^{\times}$ on all points of $\wt C$ is just $\{\pm 1\}$. 

Represent the irreducible component $\wt C=(b_p C_0, b U/U)$ with $b_p\in \GL_2(\QQ_p)$ and $b\in B_{\BA_f}^\times$. 
By translation, the common stabilizer of $D^{\times}$ on all points of $b_p C_0$ is  
$b_p \QQ_p^\times  (1+p M_2(\ZZ_p)) b_p^{-1}$. 
Assume that $\gamma\in D^\times$ fixes every point of $\wt C$. 
By $\gamma \wt C=(\gamma b_p C_0, \gamma b U/U)$, we have
$$
\gamma \in b_p \QQ_p^\times  (1+p M_2(\ZZ_p)) b_p^{-1}, \quad
\gamma b U= b U.
$$
Note that $\gamma$ acts as $v_p(q(\gamma))$ on $D_p^\times/U_p\cong \ZZ$, so 
$\gamma b U= b U$ at the place $p$ implies  $v_p(q(\gamma))=1$. 
Then the condition becomes
$$
\gamma \in \ZZ_p^\times  (1+p b_p M_2(\ZZ_p) b_p^{-1}), \quad
\gamma \in  b^p U^p (b^p)^{-1}.
$$
For convenience, we write it as 
$$
\gamma \in D^\times\cap V, \quad V=V_pV^p \subset D_{\BA_f}^\times, 
\quad V_p=  \ZZ_p^\times  (1+p b_p M_2(\ZZ_p) b_p^{-1}), \quad
V^p=  b^p U^p (b^p)^{-1}.
$$
It suffices to prove that these two conditions imply $\gamma=\pm1$. 
As they imply $q(\gamma)=\pm1$, it suffices to prove $\gamma\in \QQ^\times$. 
The proof is very similar to those in \cite[\S4.6]{Yua} and 
\cite[Prop. 4.1]{YZ}, but we reproduce it here for convenience of readers.

Assume the contrary that $\gamma\notin \QQ^\times$. Then $K=\QQ(\gamma)=\QQ+\QQ\gamma$ is a quadratic imaginary extension of $\QQ$ contained in $D$. 
Moreover, $\gamma$ lies in $V\cap K_{\BA_f}^\times$, which is an open and compact subgroup of $K_{\BA_f}^\times$.
Note that $\wh O_{K}^\times$ is the unique maximal open compact subgroup of $K_{\BA_f}^\times$.  
It follows that $\gamma\in \wh O_{K}^\times$, and thus $\gamma\in O_{K}^\times$ is a unit. 
Then it must be a root of unity as $K$ is quadratic imaginary. 
If $K$ is not isomorphic to $\QQ(\sqrt{-1})$ or $\QQ(\sqrt{-3})$, we have 
$O_K^\times =\{\pm1\}$, and the proof is finished. 

To exclude these two exceptional cases, note that 
 $\gamma$ lies in
$$
 K_p^\times \cap \big(\ZZ_p^\times  (1+p b_p M_2(\ZZ_p) b_p^{-1})\big)
= \ZZ_p^\times\cdot K_p^\times \cap \big(   1+p b_p M_2(\ZZ_p) b_p^{-1}\big)
= \ZZ_p^\times\cdot  
  \big(1+p (K_p \cap  b_p M_2(\ZZ_p)b_p^{-1}) \big).
  $$
Here we write $K_p=K\otimes_\QQ \QQ_p$ and $O_{K_p}=O_K\otimes_\ZZ \ZZ_p$. 
The intersection $K_p \cap  b_p M_2(\ZZ_p)b_p^{-1}$ is an open and compact subring of $K_p$, so it is contained in the unique maximal compact subring $O_{K_p}$ of $K_p$. It follows that 
$\gamma$ lies in $ \ZZ_p^\times (1+p O_{K_p})$. 
Then $\ZZ_p+\ZZ_p\gamma\subset \ZZ_p+ p O_{K_p}$, where the latter is an order in $O_{K_p}$. 

If $\gamma\neq \pm1$, then $\gamma$ is a root of unity of order $3,4,6$, and  $O_K=\ZZ[\gamma]=\ZZ+\ZZ\gamma$ in all these cases.
This implies $O_{K_p}=\ZZ_p+\ZZ_p \gamma$, contradicting the above result that
$\ZZ_p+\ZZ_p\gamma\subset \ZZ_p+ p O_{K_p}$. 
This finishes the proof.
\end{proof}

By the lemma, we further see that the pull-back of every connected component of $\CL_U(\CC)$ to $\CH$ is canonically isomorphic to $\Omega_{\CH/\CC}$. 
Then it is reasonable to define the \emph{Petersson metric} $\|\cdot\|_\Pet$ of $\CL_U$ by 
$$\|d\tau\|_{\mathrm{Pet}}=2\, \Im(\tau),$$
where $\tau$ is the standard coordinate function on $\CH\subset \CC$. 
The isomorphism in Lemma \ref{compatibility} is an isometry.

\begin{thm}\label{height}
Assume that $B$ is a division algebra and that 
$U=\wh O_B^\times$.
Then the normalized arithmetic intersection numbers satisfy
$$
\frac{\widehat\deg(\hat c_1(\CL_U, \|\cdot\|_\Pet)^2)}{2\deg(\CL_{U,\QQ})}
=\frac{\widehat\deg(\hat c_1(\omega_{\CX_U/\ZZ}, \|\cdot\|_\Pet)^2)}{2\deg(\omega_{\CX_{U,\QQ}/\QQ})}
=\frac{\widehat\deg(\hat c_1(\underline\omega_{\CA}, \|\cdot\|_\Fal)^2)}{2\deg(\underline\omega_{\CA,\QQ})}
+\frac12 \log d_B.
$$
\end{thm}

\begin{proof}
The first equality follows from Lemma \ref{compatibility}(2), where the isomorphism is an isometry under the Petersson metrics. 
The second equality follows from Theorem \ref{KS}.
\end{proof}

\begin{remark}
Note that \cite[Theorem 1.1]{Yua} computes the first term of the theorem (as a special case);
\cite[Theorem 1.0.5]{KRY2} computes the numerator $\widehat\deg(\hat c_1(\underline\omega_{\CA}, \|\cdot\|_\Fal)^2)$ in the third term of the theorem.
The theorem asserts that these two formulas are compatible. 
\end{remark}

\section{Image over the integral model} \label{sec image}

The goal of this section is to prove the first statement of Theorem \ref{KS}.
The key is an isomorphism from deformation theory. 

\subsection{Isomorphism from deformation theory}

Here we introduce an isomorphism related to the Kodaira--Spencer map, which comes from deformation theory. It will be used in the proof of Theorem \ref{KS}.
 
Denote by $\CX_U^{\rm sm}$ the smooth locus of $\CX_U$ over $\ZZn$, which is the maximal open substack of $\CX_U$ that is smooth over $\ZZn$.
Note that $\CX_U$ has semistable reduction over $\ZZn$, so the non-smooth locus 
$\CX_U^{\rm sing}=\CX_U\setminus \CX_U^{\rm sm}$ is 0-dimensional in a suitable sense.

Recall that $\CA\to \CX_U$ is the universal abelian scheme.
The Kodaira--Spencer map 
$$
\phi_1:\underline\Omega_{\CA} \lra \underline\Omega_{\CA^t}^\vee\otimes \omega_{\CX_U/\ZZn}.
$$
induces a morphism 
$$
\phi_3:\omega_{\CX_U/\ZZn}^\vee \lra 
\Lie(\CA^t)\otimes\Lie(\CA)$$
and a morphism
$$
\phi_4:T_{\CX_U^\sm/\ZZn} \lra \CH om( \Lie(\CA^t)^\vee, 
\Lie(\CA)  )|_{\CX_U^\sm}.
$$
By deformation theory, we know the image of this morphism. 

\begin{thm} \label{deformation}
The morphism $\phi_4$ induces an isomorphism 
$$
\phi_5:T_{\CX_U^\sm/\ZZn} \lra \CH om_{O_B}( \Lie(\CA^t)^\vee, 
\Lie(\CA)  )|_{\CX_U^\sm}.
$$
\end{thm}
\begin{proof}
The Kodaira-Spencer map $\phi_1$ induces a morphism 
$$
\phi_6:\underline\Omega_{\CA} \otimes \underline\Omega_{\CA^t}\lra \omega_{\CX_U/\ZZn}.
$$
The key point is that deformation theory implies that $\phi_6$ induces an isomorphism
$$
\phi_7:(\underline\Omega_{\CA} \otimes \underline\Omega_{\CA^t})|_{\CX_U^
\sm}/\CR\lra \omega_{\CX_U^\sm/\ZZn},
$$
where $\CR$ is the subsheaf of 
$(\underline\Omega_{\CA} \otimes \underline\Omega_{\CA^t}^\vee)|_{\CX_U^\sm}$ locally generated by 
$$
(i(\beta)^* u)\otimes v-u\otimes ((i(\beta)^t)^* v), \qquad u\in \underline\Omega_{\CA},\ v\in \underline\Omega_{\CA^t},\ \beta\in O_B.
$$
Here $i(\beta)^t:\CA^t\to \CA^t$ is the dual of $i(\beta):\CA\to \CA$ induced by the pull-back map of line bundles. 
For a serious proof of this fact, we refer to \cite[Proposition 2.3.5.2]{Lan}. 
Note that our moduli structure does not involve a polarization, so there is no polarization in the deformation, either.

To convert the isomorphism to the form in the theorem, we claim that the dual of $\phi_7$ gives an isomorphism 
$$
\phi_8:T_{\CX_U^\sm/\ZZn} \lra 
(\Lie(\CA) \otimes \Lie(\CA^t))|_{\CX_U^\sm}^{O_B},
$$
where the right-hand side denotes the subsheaf of $(\Lie(\CA) \otimes \Lie(\CA^t))|_{\CX_U^\sm}$ over $\CX_U^\sm$ locally generated by 
local sections $\sum_j x_j\otimes y_j$ satisfying 
$$
\sum_j (i(\beta) x_j)\otimes y_j=\sum_j x_j\otimes (i(\beta)^ty_j), \quad \forall\beta\in O_B.
$$
In fact, since $\omega_{\CX_U^\sm/\ZZn}$ is locally free, the exact sequence 
$$
0\lra \CR\lra (\underline\Omega_{\CA} \otimes \underline\Omega_{\CA^t})|_{\CX_U^
\sm}\lra \omega_{\CX_U^\sm/\ZZn} \lra 0
$$
induces an exact sequence
$$
0\lra T_{\CX_U^\sm/\ZZn}\lra  
(\Lie(\CA) \otimes \Lie(\CA^t))|_{\CX_U^\sm}
\lra \CR^\vee\lra  0.
$$
The remaining part can be checked by passing to an \'etale covering, and then checked Zariski locally. 

Finally, the theorem follows from the canonical isomorphism
$$
(\Lie(\CA) \otimes \Lie(\CA^t))|_{\CX_U^\sm}^{O_B}
\lra \CH om_{O_B}( \Lie(\CA^t)^\vee, 
\Lie(\CA)  )|_{\CX_U^\sm}.
$$
\end{proof}

\subsection{Image over the integral model}

Now we are ready to prove the first statement of Theorem \ref{KS}, i.e. the canonical map $\psi_3:\underline\omega_{\CA}^{\otimes 2}\to \omega_{\CX_U/\ZZn}^{\otimes 2}$ induces an isomorphism
$\underline\omega_{\CA}^{\otimes 2}\to d_B\,\omega_{\CX_U/\ZZn}^{\otimes 2}$. 
Since $\CX_U$ is regular and $\CX_U^\sing$ has codimension 2 in $\CX_U$, it suffices to prove that 
the map induces an isomorphism $\underline\omega_{\CA}^{\otimes 2}\to d_B\,\omega_{\CX_U/\ZZn}^{\otimes 2}$ over $\CX_U^\sm$. 

For simplicity, denote 
$$
\CT=\Lie(\CA), \qquad \CT'=\Lie(\CA^t)^\vee,\qquad \CN=\CH om_{O_B}( \CT', \CT ).
$$
These are coherent sheaves over $\CX_U$.
Recall that
Theorem \ref{deformation} gives an isomorphism 
$$
\phi_{5}:T_{\CX_U^\sm/\ZZn} \lra \CN|_{\CX_U^\sm}.
$$ 
A simple consequence is that $\CN|_{\CX_U^\sm}$ is a line bundle over $\CX_U^\sm$, which can also be seen by our local argument later.
There is a canonical composition
$$
\CN\otimes \CT' 
\lra \CH om_{\CO_{\CX_U}}( \CT', \CT ) \otimes \CT'
\lra \CT.
$$
Taking determinants of the composition, we obtain a morphism
$$
\CN^{\otimes 2}\otimes \det(\CT') 
\lra \det(\CT).
$$
over $\CX_U^\sm$. 
It suffices to prove that this map induces an isomorphism 
$$
\CN^{\otimes 2}\otimes \det(\CT') 
\lra d_B\det(\CT).
$$
over $\CX_U^\sm$.

It suffices to prove the base change of the result from $\ZZn$ to $\ZZ_p$ for every prime $p\nmid n$. In the following, we treat the cases $p\nmid d_B$ and $p\mid d_B$ separately.

\subsubsection*{Good prime}

We first treat case $p\nmid d_B$.
For convenience, we use subscript $p$ to indicate the base change $\otimes_\ZZ\ZZ_p$. 
Then $O_{B,p}=O_B\otimes_\ZZ \ZZ_p$
is isomorphic to $M_2(\ZZ_p)$.
Fix such an isomorphism. Take idempotents $e_1=\matrixx{1}{0}{0}{0}$, $e_2=\matrixx{0}{0}{0}{1}$ and the involution $\epsilon=\matrixx{0}{1}{1}{0}$ in $O_{B,p}$. 
The idempotents give decompositions
$$
\CT_p=e_1 \CT_p\oplus e_2\CT_p, \qquad
\CT'_p=e_1 \CT'_p\oplus e_2 \CT'_p.
$$
The involution $\epsilon$ gives isomorphisms $e_1 \CT_p\to e_2\CT_p$
and $e_1 \CT_p'\to e_2\CT_p'$.
These are locally free over $\CX_{U,p}$ since they are direct summations of locally free sheaves, and  they are line bundles by the determinant condition.

For $j=1,2$, the canonical maps
$$
\CN_p\lra \CH om_{O_B}( \CT'_p, \CT_p )\lra 
\CH om_{\CO_{\CX_{U,p}}}( e_j\CT'_p, e_j\CT_p)
\lra  (e_j\CT'_p)^\vee \otimes (e_j\CT_p)
$$
are isomorphisms. 
Then we obtain isomorphisms 
$$
\CN_p^{\otimes 2}
\lra (e_1\CT'_p)^\vee \otimes (e_1\CT_p) \otimes (e_2\CT'_p)^\vee \otimes (e_2\CT_p)
\lra \det(\CT'_p)^\vee \otimes \det(\CT_p).
$$
This proves the case $p\nmid d_B$.

\subsubsection*{Bad prime}

Here we prove the first statement of Theorem \ref{KS} at bad primes $p\mid d_B$.
Assume that $p\mid d_B$ in the following. 

To split $B$, we will pass to $\QQ_{p^2}$, the unique unramified quadratic extension of $\QQ_p$. Denote by $\ZZ_{p^2}$ the valuation ring of $\QQ_{p^2}$. 
It suffices to prove that the determinant process induces a canonical isomorphism 
$$
\CN^{\otimes 2}\otimes \det(\CT') 
\lra p\det(\CT)
$$
over the base change $\CX_U^\sm\times_{\ZZ}\ZZ_{p^2}$.
By taking an \'etale covering, we can replace $\CX_U^\sm\times_{\ZZ}\ZZ_{p^2}$ by a smooth scheme $S$ over $\ZZ_{p^2}$.

The problem is local on $S$, so it suffices to prove the corresponding result at every closed point $s\in S$ above $p$.
Denote by $R=\CO_{S,s}$ the local ring. 
Denote the $R$-modules
$$
T=\CT\times_{\CO_S}R=\Lie(\CA_R/R),\qquad
T'=\CT'\times_{\CO_S}R=\Lie(\CA_R^t/R)^\vee
$$
and
$$
N=\CN\otimes_{\CO_S} R=\Hom_{O_B}(T',T).
$$
We need to prove that the determinant process induces a canonical isomorphism
$$
N^{\otimes 2}\otimes_R \det(T') 
\lra p\det(T).
$$

For convenience, we start with some terminologies in basic algebra and basic representation.
Denote 
$$D=B\otimes_\QQ \QQ_{p},\qquad O_D=O_B\otimes_\ZZ \ZZ_{p}$$ 
Then $O_D$ is the unique maximal order of the division quaternion algebra $D$.  

By a \emph{left} (resp. \emph{right}) \emph{$O_D$-module of rank 2 over $R$}, we mean a pair $(M,i)$,  where $M$ is a free $R$-module of rank 2 and $i:O_D\to \End_R(M)$ is a 
$\ZZ_{p}$-linear homomorphism (resp. anti-homomorphism).
A \emph{one-sided} $O_D$-module of rank 2 over $R$ means either a left or a right $O_D$-module of rank 2 over $R$;
\emph{an $O_D$-module of rank 2 over $R$} means a left $O_D$-module of rank 2 over $R$.
In either case, the determinant condition means $\det(i(\beta))=q(\beta)$ for any $\beta\in O_D$.

Now we introduce two operations which keep the determinant condition.
The first is the conjugation action. 
Let $\gamma\in D$ be a nonzero element, which gives the conjugation action
$$\ad(\gamma):O_D\lra O_D,\qquad
\beta\longmapsto \gamma\beta\gamma^{-1}.$$
If $(M,i)$ is a one-sided $O_D$-module of rank 2 over $R$, then 
so is $(M,i\circ \ad(\gamma))$.

The second is the dual module. 
If $(M,i)$ is a one-sided $O_D$-module of rank 2 over $R$, 
define the \emph{dual} $(M^\vee,i^\vee)$ as follows. 
The underlying $R$-module $M^\vee=\Hom_R(M,R)$, and the action of $\beta\in O_D$ sends $\ell\in \Hom_R(M,R)$ to $\ell\circ i(\beta)$. 
Then $(M^\vee,i^\vee)$ is also a one-sided $O_D$-module of rank 2 over $R$.
The dual process switch left modules and right modules. 
   
As a basic result, if $(M,i)$ is a one-sided $O_D$-module of rank 2 over $R$ satisfying the determinant condition, 
then $(M^\vee,i^\vee)$ is isomorphic to $(M, i\circ \iota)$
as a one-sided $O_D$-module. Here $\iota:O_D\to O_D$ is the main involution. 
In fact, consider the $R$-bilinear pairing
$$
M\times M\longrightarrow R, \quad (x,y)\longmapsto \det(x,y). 
$$
Here we fix an isomorphic $M\simeq R^2$ and view elements of $M$ as column vectors, so $\det(x,y)$ is the determinant of a $2\times 2$ matrix.
For any nonzero $\beta\in D$, the determinant condition givesÊ
$$\det(i(\beta)x, y)
=q(\beta)\det(i(\beta)^{-1}i(\beta)x, i(\beta)^{-1}y)
=\det(x, q(\beta)i(\beta)^{-1}y)
=\det(x, i(\beta^\iota)y).
$$
The pairing induces an isomorphism $M^\vee\to M$, under which $i^\vee$ becomes $i\circ \iota$.

Return to
$$
T=\Lie(\CA_R/R),\qquad
T'=\Lie(\CA_R^t/R)^\vee
$$
Then $(T,i)$ is a left $O_D$-module of rank 2 over $R$ satisfying the determinant condition, where $i:O_D\to \End_R(T)$ is as in the moduli structure.
We claim that $T'$ with the natural action is a left $O_D$-module of rank 2 over $R$ isomorphic to $(T,i\circ \ad(\mu))$, where $\mu\in O_B$ is the element used to define the positive involution $\beta^*=\mu^{-1} \beta^\iota \mu=\mu \beta^\iota \mu^{-1}$ on $B$. 

In fact,  the principal polarization induces an isomorphism $T\to \Lie(\CA_R^t/R)$.
By compatibility with the Rosati involution, $\Lie(\CA_R^t/R)$ with the dual action of $O_D$ is isomorphic to $(T,i\circ *)$.
Then $T'=\Lie(\CA_R^t/R)^\vee$ is isomorphic to $(T,i\circ *\circ \iota)=(T,i\circ \ad(\mu))$.

Now the theorem follows from the following basic result on abstract 
 $O_D$-modules of rank 2 over $R$ satisfying the determinant condition.

\begin{lem}
Let $R$ be a flat noetherian local ring over $\ZZ_{p^2}$ such that $pR$ is a prime ideal of $R$. 
Let $(T,i)$ be a left $O_D$-module of rank 2 over $R$ satisfying the determinant condition.
Let $(T',i')=(T,i\circ \ad(\mu))$ be another left $O_D$-module over $R$, where $\mu\in O_D$ is any element such that $\ord_p(q(\mu))=1$.
Then the following hold. 
\begin{enumerate}[(1)]
\item
The $R$-module $N=\Hom_{O_B}(T',T)$ is free or rank 1.
\item
Let 
$$
\psi:N^{\otimes 2}\otimes_R \det(T') \lra \det(T)
$$
be the determinant of the composition map
$$
N \otimes_R T' \lra \Hom_{R}(T',T) \otimes_R T'  \lra T.
$$
Then $\psi$ is injective with
$\Im(\psi)=p\det(T)$.
\end{enumerate}
\end{lem}
\begin{proof}
First, we claim that there is an $R$-linear ring isomorphism
$$
\tau:O_D\otimes_{\ZZ_p} R\lra \matrixx{R}{R}{pR}{R}.
$$
For this, it suffices to check the case $R=\ZZ_{p^2}$. 
Write $O_D=\ZZ_{p^2}+\ZZ_{p^2} j$ with $j\in O_D$ satisfying
$$j^\iota=-j, \quad j^2=p,\quad j x= \bar x j,\ \forall x\in \ZZ_{p^2}.$$ 
Here $\bar x$ denotes the image of $x$ under the nontrivial element of $\Gal(\QQ_{p^2}/\QQ_{p})$, which is compatible with the main involution of $D$.
Take an injective homomorphism 
$$
\rho_1:O_D\lra M_2(\ZZ_{p^2}),\quad x\longmapsto \matrixx{x}{}{}{\bar x},
\quad j\longmapsto \matrixx{}{1}{p}{},\quad x\in \ZZ_{p^2}.
$$
It induces an injective homomorphism
$$
\rho_2:O_D\otimes_{\ZZ_p} \ZZ_{p^2} \lra M_2(\ZZ_{p^2}),\quad
x\otimes a\longmapsto \matrixx{ax}{}{}{a\bar x},
\quad j\longmapsto \matrixx{}{1}{p}{},\quad x\otimes a\in \ZZ_{p^2}\otimes_{\ZZ_p} \ZZ_{p^2}.
$$
As $\ZZ_{p^2}$ is unramified over $\ZZ_p$, 
there is an isomorphism 
$$
 \ZZ_{p^2}\otimes_{\ZZ_p}  \ZZ_{p^2}\lra  \ZZ_{p^2}\oplus \ZZ_{p^2}, \quad
 x\otimes a\longmapsto (ax,a\bar x).
$$
Then we see that the image of $\rho_2$ is as required.

Second, we claim that up to isomorphism, there are exactly two 
$O_D$-modules $M$ of rank 2 over $R$ satisfying the determinant condition.
Identify $O_D\otimes_{\ZZ_p}R$ with the matrix ring via $\tau$. 
Then $O_D\otimes_{\ZZ_p}R$ contains elements
$$e_1=\matrixx{1}{0}{0}{0},\quad
e_2=\matrixx{0}{0}{0}{1}, \quad
j=\matrixx{0}{1}{p}{0}.
$$ 
We can write
\begin{equation*} \label{aaaa}
M=e_1M \oplus e_2 M, \quad e_1M=Rx, \quad e_2M=Ry
\end{equation*}
as $R$-modules. Note that $e_1M, e_2M$ are projective and thus free over $R$, and they have rank $1$ by the determinant condition. The elements $x,y$ are their bases.
As $je_1=e_2j$ and $je_2=e_1j$, the action of $j$ on $M$ switches $e_1M, e_2M$.
Write the action
$$jx=ay,\quad  jy=bx, \quad a,b\in R.$$ 
Then we have $ab=p$ by $j^2=p$.
By assumption, $pR$ is a prime ideal, so either $a\in R^\times$ or $b\in R^\times$. 
If $b\in R^\times$, replacing $x$ by $b^{-1}x$ in the above, we obtain
\begin{equation*} \label{bbbb}
jx=py,\quad  jy=x.
\end{equation*}
This is the standard representation.
If $a\in R^\times$, replacing $y$ by $a^{-1}y$ in the above, we obtain
$$jx=y,\quad  jy=px.$$ 
These give the two $O_D$-modules, which are not isomorphic by considering surjectivity of $j:e_1M\to e_2M$. 

Third, let $T$ and $T'$ be the $O_D$-modules as in the lemma. 
Then up to isomorphism, $T$ and $T'$ are exactly the two distinct $O_D$-modules in the above classification.
In fact, as $O_D$ is the maximal order and $\ord_p q(j)=\ord_p q(\mu)$, the quotient
$j\mu^{-1}$ lies in $O_D^\times$.
Then in the definition of $T'$, we can assume $\mu=j$. 
As the action $\ad(j):O_B\to O_B$ switches $e_1, e_2$, we see that composing $O_B\to \End(M)$ with $\ad(j)$ switches the two representations above. 
This proves that $T$ and $T'$ are the two distinct $O_D$-modules.

Fourth, we can finally prove the lemma. 
We only consider the case that $T$ is the above standard representation $M$, since the other case is similar. 
Then $T$ has  a presentation
\begin{equation*} 
T=e_1T \oplus e_2 T, \quad e_1T=Rx, \quad e_2T=Ry,\quad 
jx=py,\quad  jy=x.
\end{equation*}
Then $T'$ has a presentation
\begin{equation*} 
T'=e_1T' \oplus e_2 T', \quad e_1T'=Rx', \quad e_2T'=Ry',\quad 
jx'=y',\quad  jy'=px'.
\end{equation*}
Any $f\in \Hom_{O_B}(T',T)$ sends $e_iT'$ to $e_iT$ for $i=1,2$.
Write $f(x')=ax$ and $f(y')=by$ with $a,b\in R$. 
The action of $j$ gives $f(y')=apy$ and $f(px')=bx$.
Comparing these four relations gives $b=pa$. 
One can check that any pair $(a, pa)$ with $a\in R$ gives a unique homomorphism $f$.
Then the canonical maps
$$\Hom_{O_B}(T',T)\lra  \Hom_{R}(e_1T',e_1T),$$
$$\Hom_{O_B}(T',T)\lra p\,\Hom_{R}(e_2T',e_2T)$$
are isomorphisms.
Hence, the natural map
$$\Hom_{O_B}(T',T)\otimes_R T' \lra T$$
induces an isomorphism
$$
\Hom_{O_B}(T',T)\otimes_R T' \lra e_1T\oplus pe_2T.
$$
Taking determinants, we have canonical isomorphisms
$$
(\Hom_{O_B}(T',T))^{\otimes 2}\otimes_R \det(T') \lra (e_1T)\otimes_R (pe_2T)
\lra p\det(T).
$$
This finishes the proof.
\end{proof}

\section{Comparison of the metrics}  \label{sec metric}

The goal of this section is to prove the compatibility of metrics in Theorem \ref{KS}. 
Then the situation is complex analytic, and can be performed over the universal abelian surface over the upper half plane.
In \S\ref{sec tangent}, we use \v Cech cohomology to describe the map between tangent spaces induced by a polarization, which will be used in the proof afterwards.
In \S\ref{sec explicit}, we compute the Kodaira--Spencer map over the upper half plane explicitly.
Then the comparison of the metrics is easily done in \S\ref{sec comparison}.

\subsection{Maps between complex tangent spaces} \label{sec tangent}

Let $A=V/\Lambda$ be a $g$-dimensional complex abelian variety. 
Here $V$ is a complex vector space of dimension $g$, and $\Lambda$ is a full lattice in $V$. We have identification $V=\Lambda\otimes_\ZZ\RR$ and canonical isomorphism $\Lie(A)\simeq V$. 
Let $\lambda:A\to A^t$ be a polarization induced by a Riemann form 
$$E:\Lambda\times \Lambda\lra \ZZ.$$ 
The polarization is not assumed to be principal here.
We refer to \cite[\S1-3,\S9]{Mum} for basics of complex abelian varieties.

This induces isomorphisms
$$
V\lra \Lie(A) \lra \Lie(A^t) \lra H^1(A, \CO_A),
$$
where the second arrow is by the polarization and the third arrow is the canonical isomorphism by deformation of line bundles over $A$. The goal here is to describe this composition explicitly. 

\subsection*{\v Cech cohomology}

The Hodge theory gives a split exact sequence 
$$
0\lra H^0(A, \Omega_{A/\CC}) \lra H^1(A, \CC)
\lra H^1(A, \CO_{A})
\lra 0.
$$
We have canonical isomorphisms
$$
H^1(A, \CC)\simeq \Hom_{\ZZ}(H_1(A, \ZZ),\CC)
\simeq \Hom_{\ZZ}(\Lambda,\CC)
\simeq \Hom_{\RR}(V,\CC)
\simeq \Hom_{\CC}(V\otimes_\RR\CC,\CC),
$$
together with a canonical decomposition
$$
\Hom_{\CC}(V\otimes_\RR\CC,\CC)
\simeq \Hom_{\CC\text{-linear}}(V,\CC) \oplus  \Hom_{\CC\text{-semilinear}}(V,\CC).
$$
This decomposition corresponds to canonical isomorphisms
$$
H^0(A, \Omega_{A/\CC})
\simeq \Hom_{\CC\text{-linear}}(V,\CC)
$$
and 
$$
H^1(A, \CO_{A})
\simeq \Hom_{\CC\text{-semilinear}}(V,\CC).
$$
Although these maps are canonical, they are not explicit or convenient enough for our purpose, so we will not use them in the above setting. 

Alternatively, we will use \v Cech cohomology for $H^1(A,\CO_A)$ to introduce an explicit canonical map 
$$
h:\Hom_{\ZZ}(\Lambda,\CC) \lra H^1(A,\CO_A).
$$  
More generally, let $\CF$ be either the sheaf $\CO$ or the sheaf $\CO^\times$ in the complex analytic setting. 
We are going to introduce an explicit homomorphism
$$
\delta: H^1(\Lambda, \CF(V)) \lra H^1(A,\CF_A). 
$$
The homomorphism is an isomorphism for both cases, but we do not need this fact.
Our approach fits the general Hochschild--Serre spectral sequence for the Galois cover $V\to A$, but we do not need the spectral sequence, either.

We first need an open cover of $A$ coming from the universal cover $V$. 
Take a set $\{U_t\}_{t\in I}$ of open subsets $U_t$ of $V$ satisfying the following conditions:
\begin{itemize}
\item[(1)] every composition $U_t\to V\to A$ is injective, 
\item[(2)] the induced map $\cup_{t\in I} U_t\to A$ is surjective, 
\item[(3)] for any $t,t'\in I$ (allowing $t=t'$), the difference $U_{t'}-U_{t}=\{z'-z:z\in U_{t}, z'\in U_{t'}\}$, as a subset of $V$, contains at most one point of $\Lambda$. 
This point is denoted by $c_{t,t'}\in \Lambda$ if it exists. 
\end{itemize}
Denote $\bar U_t$ the image of $U_t\to A$. 
Then $\{\bar U_t\}_{t\in I}$ forms an open cover of $A$.
We call it an \emph{admissible cover} of $A$. 
We are going to use this open cover to compute \v Cech cohomology.  

Now we define the explicit map 
$$
\delta: H^1(\Lambda, \CF(V)) \lra H^1(A,\CF_A). 
$$ 
For a cross-homomorphism $\alpha:\Lambda\to \CF(V)$, we define a \v Cech cocycle $\delta(\alpha)$ whose component $\delta(\alpha)_{t,t'}\in \CF(\bar U_{t,t'})$ on the overlap $\bar U_{t,t'}=\bar U_{t}\cap \bar U_{t'}$ is defined by the image of $\alpha(c_{t,t'})$ under the composition
$$\CF(V)\lra \CF(U_t)\lra \CF(\bar U_t)\lra \CF(\bar U_{t,t'}).$$ 
One checks that this is a \v Cech cocycle. 
Then we set 
$$
\delta(\alpha):=(\delta(\alpha)_{t,t'}) \in H^1(A,\CF_A).
$$

Now we define the map
$$
h:\Hom_{\ZZ}(\Lambda,\CC) \lra H^1(A,\CO_A)
$$
as the composition
$$
\Hom_{\ZZ}(\Lambda,\CC)\lra H^1(\Lambda,\CC) \lra H^1(\Lambda,\CO(V)) \stackrel{\delta}{\lra} H^1(A,\CO_A).
$$ 
Here the first arrow is the canonical isomorphism by the fact that $\Lambda$ acts trivially on $\CC$ (viewed as the space of constant functions), and the second arrow is the map induced by the natural map $\CC\to \CO(V)$.
Similarly, we define the map 
$$
h^\times:\Hom_{\ZZ}(\Lambda,\CC^\times) \lra H^1(A,\CO_A^\times).
$$ 
as the composition
$$
\Hom_{\ZZ}(\Lambda,\CC^\times)\lra H^1(\Lambda,\CC^\times) \lra H^1(\Lambda,\CO^\times(V)) \stackrel{\delta}{\lra} H^1(A,\CO_A^\times).
$$

\subsubsection*{Maps between tangent spaces}

Recall that the polarization $\lambda:A\to A^t$ is induced by the Riemann form $E:\Lambda\times\Lambda\to \ZZ$. Now we can present our main result of this subsection.

\begin{lem}\label{tangent space}
The composition
$$V\lra \Lie(A)\stackrel{\lambda}{\lra} \Lie(A^t) \lra H^1(A,\CO_A)$$
is given by 
$$
z\longmapsto 2\pi i \, h(E(z,\cdot)),
$$ 
where 
$E(z,\cdot):\Lambda\to \CC$ is viewed as an element of $\Hom_\ZZ(\Lambda, \CC)$, 
and $h: \Hom_\ZZ(\Lambda, \CC) \to H^1(A,\CO_A)$
is the canonical map defined above. 
\end{lem}

\begin{proof}

Let us first describe the polarization $\lambda: A\to  A^t$ 
in terms of the Riemann form $E:\Lambda\times\Lambda\to \ZZ$. 
Denote by $H:V\times V\to \CC$  the hermitian form corresponding to $E$.
By \cite[p. 20, \S2]{Mum}, there is a line bundle $L(H,\alpha)$ over $A$ depending on the choice of a map
$$\alpha:\Lambda\lra S^1=\{z\in \CC:|z|=1\}$$
satisfying 
$$
\alpha(u_1+u_2)=\alpha(u_1)\alpha(u_2)e^{\pi i E(u_1, u_2)}. 
$$
This gives the polarization map 
$$
\lambda: A\lra  A^t,\quad x\longmapsto T_x^*(L(H, \alpha)) \otimes L(H,\alpha)^\vee.
$$
The line bundle $L(H,\alpha)$ depends on the choice of $\alpha$, but its cohomology class is given by $H$ and thus independent of $\alpha$. As a consequence, 
the polarization $\lambda$ is also  independent of $\alpha$.

By \cite[\S9, p. 84]{Mum}, we actually have 
$$
\lambda([z])=L(0,\gamma_z),\quad z\in V.
$$
Here $[z]\in A$ denotes the image of $z$ under $V\to A$, and  
$$\gamma_z:\Lambda\lra S^1,\quad
\gamma_z(u)=e^{2\pi i E(z,u)}.$$
Note that $\gamma_z$ lies in $\Hom(\Lambda, \CC^\times)$.
By the loc. cit., the class of the line bundle $L(0,\gamma_z)$ in 
$\Pic(A)\simeq H^1(A, \CO_A^\times)$ is exactly $h^\times(\gamma_z)$ via 
$$
h^\times: \Hom(\Lambda,\CC^\times)\lra H^1(A, \CO_A^\times). 
$$
In summary, the natural composition 
$$A\stackrel{\lambda}{\lra} A^t\lra  \Pic^0(A)\lra  \Pic(A)\lra  H^1(A,\CO_A^\times)$$ 
is given by
$$
\lambda_1: A\lra H^1(A,\CO_A^\times), \quad [z]\longmapsto h^\times(\gamma_z).
$$

We can obtain the map 
$$
\Lie(A)\lra H^1(A,\CO_A)
$$
by differentiating the map 
$$
\lambda_1: A\lra H^1(A,\CO_A^\times).
$$
This corresponds to deformations of line bundles over $A$ in algebraic geometry. 
For that, fix $z\in V$, and take $\epsilon$ to be a real number converging to 0. 
Then the above argument gives 
$$\lambda_1([\epsilon z])=h^\times(\gamma_{\epsilon z})=
h^\times(e^{2\pi i E(\epsilon z,\cdot)})
=h^\times(1+2\pi i E(z,\cdot) \epsilon+O(\epsilon^2)).$$
The linear term in $\epsilon$ gives 
$$
(d\lambda_1)(z) = 2\pi i h(E(z,\cdot)).
$$
This finishes the proof.
\end{proof}

\subsection{Explicit map in the complex setting} \label{sec explicit}

Consider the universal abelian surface $\pi:\CA\to \CH$
given by 
$$\CA=O_B\backslash (\CH\times \CC^2),$$ 
where 
$O_B$ acts on $\CH\times \CC^2$ by 
$$\beta (\tau, z)  =(\tau, z+\sigma(\beta) (\tau,1)^t ).$$
Here $(\tau,1)^t\in \CC^2$ is viewed as a column vector, and $\sigma(\beta) (\tau,1)^t$ is the matrix multiplication via the identification $\sigma:B_\RR\to M_2(\RR)$.
For simplicity, we will write
$\beta (\tau,1)^t$ for $\sigma(\beta) (\tau,1)^t$.

For each $\tau\in \CH$, we have a canonical uniformization
$$\CA_\tau= \CC^2/\Lambda_\tau,\quad \Lambda_\tau=O_B (\tau,1)^t.$$
Note that the map $O_B\to O_B (\tau,1)^t=\Lambda_\tau$ is injective, so $\Lambda_\tau$ has rank 4 over $\ZZ$.  
The space $\Lambda_\tau\otimes_\ZZ\RR=M_2(\RR)(\tau,1)^t=\CC^2$. 
So $\Lambda_\tau$ is indeed a full lattice of $\CC^2$.

The complex torus $\CA_\tau$ has a canonical positive Riemann form as follows. 
Recall that we have an element $\mu\in O_B$ such that $\mu^2=-d_B$. 
It induces a symplectic pairing
$$
E: \Lambda_\tau \times \Lambda_\tau \lra \ZZ, \qquad 
E(\beta (\tau,1)^t,\beta' (\tau,1)^t)=-\tr(\mu^{-1} \beta   \beta'^\iota).
$$
By \cite[Lemma 43.6.16, Lemma 43.6.22]{Voi}, this is indeed a self-dual positive Riemann form.
It gives principal polarizations 
$$
\lambda:\CA_\tau\lra \CA_\tau^t,\qquad
\lambda:\CA\lra \CA^t
$$

Denote by $\Lambda$ the inverse image in $\CH\times \CC^2$ of the identity section of $\CA\to \CH$.
Then we have 
$$\Lambda=\bigcup_{\tau\in \CH} \Lambda_\tau=\bigcup_{\tau\in \CH} O_B (\tau,1)^t
=\{\beta (\tau,1)^t: \beta\in O_B, \tau\in \CH\}.$$
We have a biholomorphic map
$$\Lambda \lra \CH\times O_B,\quad \beta (\tau,1)^t\longmapsto (\tau,\beta),$$
and a natural bijection $\pi_0(\Lambda) \to O_B$.

\subsubsection*{Explicit Kodaira--Spencer map}

Now we are ready to compute the Kodaira--Spencer map
$$
\phi:\pi_* \Omega_{\CA/\CH} \lra  \Lie(\CA/\CH)\otimes \Omega_{\CH/\CC}.
$$
Similar to the algebraic setting, we denote
$$
\underline\Omega_{\CA/\CH}=\pi_* \Omega_{\CA/\CH},\qquad
\underline\omega_{\CA/\CH}=\pi_* \omega_{\CA/\CH}=\det\underline\Omega_\CA.
$$
Let $\tau$ be the standard coordinate function of $\CH\subset \CC$.
Let $z=(z_1,z_2)$ be the standard coordinate functions of $\CC^2$. 
Then we have
$$\Omega_{\CA/\CH}=\CO_{\CA}dz_1+\CO_{\CA}dz_2,\quad
\underline\Omega_{\CA/\CH}=\CO_{\CH}dz_1+\CO_{\CH}dz_2,\quad
\underline\omega_{\CA/\CH}=\CO_{\CH}dz_1\wedge dz_2,\quad
\Omega_{\CH/\CC}=\CO_\CH d\tau.
$$ 

In the following, we write 
$\mu=\matrixx{a}{b}{c}{d}$ under $\sigma:B_\RR\to M_2(\RR)$. 
The explicit Kodaira--Spencer map is as follows. 

\begin{thm} \label{analytic}
The Kodaira--Spencer map
$$
\phi:\pi_* \Omega_{\CA/\CH} \lra \Lie(\CA/\CH)\otimes \Omega_{\CH/\CC}
$$
gives
$$
\phi(dz_1)=\frac{1}{2\pi i}(b\frac{\partial}{\partial z_1}+d \frac{\partial}{\partial z_2})\otimes d\tau,
$$ 
$$
\phi(dz_2)=-\frac{1}{2\pi i}(a\frac{\partial}{\partial z_1}+c \frac{\partial}{\partial z_2})\otimes d\tau.
$$
Therefore, the map 
$$
\psi: \underline\omega_{\CA/\CH}^{\otimes 2} \lra \omega_{\CH/\CC}^{\otimes 2}
$$ 
induced by $\det(\phi)$ gives 
$$
\psi((dz_1\wedge dz_2)^{\otimes 2})= 
\frac{d_B}{(2\pi i)^2}  (d\tau)^{\otimes 2}.
$$

\end{thm}

\begin{proof}
It is easy to prove the second statement by the first statement. In fact, the map
$$
\det(\phi): \det\pi_*\Omega_{\CA/\CH} \lra \det\Lie(\CA/\CH)\otimes \omega_{\CH/\CC}^{\otimes 2}
$$
is given by 
$$
dz_1\wedge dz_2\longmapsto\frac{ad-bc}{(2\pi i)^2}(\frac{\partial}{\partial z_1}\wedge\frac{\partial}{\partial z_2})\otimes (d\tau)^{\otimes 2}.
$$
This gives the second statement  by $ad-bc=q(\mu)=-\mu^2=d_B$. 

Now we prove the first statement.
Recall that 
$$
\phi_0:\pi_* \Omega_{\CA/\CH} \stackrel{}{\lra} R^1\pi_* (\pi^* \Omega_{\CH/\CC}) = R^1\pi_* \CO_\CA\otimes \Omega_{\CH/\CC}
$$
is the connecting map induced by the short exact sequence
$$
0\lra \pi^* \Omega_{\CH/\CC} \lra \Omega_{\CA/\CC} \lra \Omega_{\CA/\CH} \lra 0.
$$
Fix a point $\tau_0\in\CH$. 
We claim that on the fibers above $\tau_0$,
$$
\phi_0(dz_j)|_{\tau_0}= h(\ell_{j1}) \otimes d\tau|_{\tau_0},\quad j=1,2.
$$
Here the canonical map
$$
h:\Hom_\ZZ(\Lambda_{\tau_0},\CC)\lra H^1(\CA_{\tau_0}, \CO_{\CA_{\tau_0}})
$$
is defined before, and
$$
\ell_{jk}: \Lambda_{\tau_0}=\sigma(O_B) (\tau_0, 1)^t \lra\sigma(O_B) \lra \RR
$$
picks up the $(j,k)$-coefficient of $\sigma(O_B)\subset M_2(\RR)$.

We will use \v Cech cohomology to prove the claim.
Similar to the case of a single abelian variety in \S\ref{sec tangent},  we can have a notion of admissible cover of $\CA$ by the universal cover $\CH\times \CC^2$. 
In fact, take a set $\{U_t\}_{t\in I}$ of open subsets $U_t$ of $\CH\times \CC^2$ satisfying the following conditions:
\begin{itemize}
\item[(1)] every composition $U_t\to  \CH\times \CC^2\to \CA$ is injective, 
\item[(2)] the induced map $\cup_tU_t\to \CA$ is surjective, 
\item[(3)] for any $t,t'\in I$, the difference $U_{t'}-U_{t}=\{(\tau,z'-z):(\tau,z)\in U_{t}, (\tau,z')\in U_{t'}\}$ in $\CH\times \CC^2$ intersects at most one connected component of $\Lambda=\cup_\tau\Lambda_\tau$. 
Denote by $\beta_{t,t'}\in O_B$ the element representing this connected component (if non-empty) via the natural bijection $\pi_0(\Lambda)\to O_B$. 
\end{itemize}
Denote by $\bar U_t$ the image of $U_t\to \CA$. 
Then $\{\bar U_t\}_{t\in I}$ forms an admissible open cover of $\CA$.

Now we are ready to use \v Cech cohomology with respect to the admissible cover to compute $\phi_0(dz_j)$,  the image of the connecting map induced by
$$
0\lra \pi^* \Omega_{\CH/\CC} \lra \Omega_{\CA/\CC} \lra \Omega_{\CA/\CH} \lra 0.
$$
Note that the section $dz_j$ of $\Omega_{\CA/\CC}$ over each $U_t$ lifts the section $dz_j$ of $\pi_* \Omega_{\CA/\CH}$.
Denote by $(dz_j)_t$ the pushward of $dz_j$ via $U_t\to \bar U_t$.
Now we compute $(dz_j)_{t'}-(dz_j)_t$ on the overlap $\bar U_{t,t'}=\bar U_t\cap \bar U_{t'}$ (if nonempty). 
This overlap induces an isomorphism $U_{t,t'}\to U_{t',t}$, where $U_{t,t'}$ (resp. $U_{t',t}$) is the subset of $U_t$ (resp. $U_t'$) bijective to $\bar U_{t,t'}$.
By definition, the isomorphism $U_{t,t'}\to U_{t',t}$ maps $(\tau, z)$ to $(\tau,z+ \beta_{t,t'}(\tau,1)^t)$. Thus the pull-back of the function $z_j\in \CO(U_{t',t})$
to $\CO(U_{t,t'})$ becomes 
$$z_j'=z_j+ \ell_{j1}(\beta_{t,t'})\,\tau+\ell_{j2}(\beta_{t,t'}),$$
where we write $\sigma(\beta_{t,t'})=\matrixx{\ell_{11}(\beta_{t,t'})}{\ell_{12}(\beta_{t,t'})}{\ell_{21}(\beta_{t,t'})}{\ell_{22}(\beta_{t,t'})}$. 
As a consequence, the pull-back of 
$(dz_j)_{t'}-(dz_j)_{t}$ to $U_{t,t'}$ is just 
$$
dz_j'-dz_j = d(z_j+ \ell_{j1}(\beta_{t,t'})\,\tau+\ell_{j2}(\beta_{t,t'}))-dz_j
= \ell_{j1}(\beta_{t,t'})\ d\tau.
$$
As a consequence, the class $\phi_0(dz_j)$ in 
$R^1\pi_* \CO_\CA\otimes\Omega_{\CH/\CC}$ is represented by the \v Cech cocycle 
$$(\ell_{j1}(\beta_{t,t'})\, d\tau)_{t,t'\in J}=(\ell_{j1}(\beta_{t,t'}))_{t,t'\in J}\otimes d\tau.$$
Hence, the restriction of $\phi_0(dz_j)$ to the fiber $\CA_{\tau_0}$ above $\tau_0\in \CH$ is exactly
$$
\phi_0(dz_j)|_{\tau_0}= h(\ell_{j1})\otimes d\tau|_{\tau_0} \in 
H^1(\CA_{\tau_0}, \CO_{\CA_{\tau_0}})\otimes \Omega_{\CH/\CC}(\tau_0).
$$
This proves the claim.

By Lemma \ref{tangent space}, the composition
$$\CC^2\simeq \Lie(\CA_{\tau_0})\lra \Lie(\CA_{\tau_0}^t) \lra H^1(\CA_{\tau_0},\CO_{\CA_{\tau_0}})$$
is given by 
$$
\rho:z\longmapsto 2\pi i \, h(E(z,\cdot)).
$$ 
We are going to find a $w_j\in \CC^2$ such that 
$$h(\ell_{j1})=(2\pi i)^{-1}\rho(w_j)=h(E(w_j,\cdot)).$$
It suffices to find a $w_j\in \CC^2$ satisfying $E(w_j,\cdot)=\ell_{j1}$ in 
$\Hom_\ZZ(\Lambda_{\tau_0},\CC)$.
Recall that the Riemann form over $\Lambda_{\tau_0}=O_B ({\tau_0},1)^t$ is given by
$$
E: \Lambda_{\tau_0} \times \Lambda_{\tau_0} \lra \ZZ, \qquad 
E(\beta ({\tau_0},1)^t,\beta' ({\tau_0},1)^t)=-\tr(\mu^{-1} \beta \beta'^\iota).
$$
Set $w_j=\beta_j ({\tau_0},1)^t$ with $\beta_j\in M_2(\RR)$.
Then the equation becomes
$$
-\tr(\mu^{-1} \beta_j \beta'^\iota)=\ell_{j1}(\beta'),\quad \forall \beta'\in M_2(\RR).
$$
Recall that $\mu=\matrixx{a}{b}{c}{d}$ under $\sigma:B_\RR\to M_2(\RR)$. 
Some matrix calculations give a solution
$$
\beta_j=
\begin{cases}
\matrixx{0}{b}{0}{d}&j=1,\\
\matrixx{0}{-a}{0}{-c} &j=2.
\end{cases}
$$
The corresponding 
$$
w_j=\beta_j(\tau_0,1)^t=
\begin{cases}
(b, d)^t&j=1,\\
(-a, -c)^t &j=2.
\end{cases}
$$
These vectors in $\CC^2$ are written respectively as $\displaystyle b\frac{\partial}{\partial z_1}+d \frac{\partial}{\partial z_2}$
and
$\displaystyle -a\frac{\partial}{\partial z_1}-c \frac{\partial}{\partial z_2}$ under the identification $\CC^2\simeq \Lie(\CA_{\tau_0})$.

Finally, the relation $h(\ell_{j1})=(2\pi i)^{-1}\rho(w_j)$ gives
$$
\phi_0(dz_1)|_{\tau_0}= \frac{1}{2\pi i}\rho(b\frac{\partial}{\partial z_1}+d \frac{\partial}{\partial z_2})\otimes d\tau\in  
H^1(\CA_{\tau_0}, \CO_{\CA_{\tau_0}})\otimes \Omega_{\CH/\CC}(\tau_0),
$$
and
$$
\phi_0(dz_2)|_{\tau_0}=\frac{1}{2\pi i}\rho(-a\frac{\partial}{\partial z_1}-c \frac{\partial}{\partial z_2})\otimes d\tau\in  
H^1(\CA_{\tau_0}, \CO_{\CA_{\tau_0}})\otimes \Omega_{\CH/\CC}(\tau_0),
$$
It follows that the map 
$$
\phi:\pi_* \Omega_{\CA/\CH} \lra \Lie(\CA/\CH)\otimes \Omega_{\CH/\CC},
$$
which satisfies $\phi_0=(\rho\otimes 1)\circ\phi$,
 gives
$$
\phi(dz_1)|_{\tau_0}
=\frac{1}{2\pi i}(b\frac{\partial}{\partial z_1}+d \frac{\partial}{\partial z_2})\otimes d\tau \in \Lie(\CA_{\tau_0})
\otimes \Omega_{\CH/\CC}(\tau_0),
$$
and 
$$
\phi(dz_2)|_{\tau_0}=-\frac{1}{2\pi i}(a\frac{\partial}{\partial z_1}+c \frac{\partial}{\partial z_2})\otimes d\tau
\in \Lie(\CA_{\tau_0})
\otimes \Omega_{\CH/\CC}(\tau_0).
$$
This finishes the proof.
\end{proof}

\subsection{Comparison of the metrics} \label{sec comparison}

Now we are ready to prove the second statement of Theorem \ref{KS}.
Every connected component of $\CX_U(\CC)$ is a quotient of $\CH$ by a discrete subgroup of $O_B^\times$. 
The pull-back abelian schemes of $\CA(\CC)\to \CX_U(\CC)$
to $\CH$ (via connected components of $\CX_U(\CC)$) are the universal abelian variety $\CA\to \CH$ in the complex setting described above.
The Hodge bundles, the canonical bundle, the Kodaira--Spencer maps, and the metrics are compatible with the pull-back.
So we only need to compare the metric in this complex setting.

The Faltings metric on $\underline\omega_{\CA/\CH}$ is given by 
$$
\| \alpha_{\tau_0} \|^2_{\Fal} =  \frac{1}{(2\pi)^2} \left| \int_{\CA_{\tau_0}} \alpha_{\tau_0}\wedge \bar\alpha_{\tau_0} \right|,
$$
where $\alpha_{\tau_0}\in \underline\omega_{\CA/\CH}(\tau_0)$ for $\tau_0\in\CH$. 
The Peterson metric on $\omega_{\CH/\CC}$ is given by 
$$
\|d\tau\|_{\rm Pet} =  2\, {\rm Im}(\tau).
$$
Now we will see that these two metrics are equal under the Kodaira--Spencer map.

As before, let $z_1,z_2$ be coordinates of $\CC^2$ as in $\CA_{\tau_0}=\CC^2/O_B(\tau_0,1)^t$. 
By definition, we have
$$
\| dz_1\wedge dz_2\|^2_{\Fal} (\tau_0)=  \frac{1}{(2\pi)^2} \left| \int_{\CA_{\tau_0}} 
dz_1\wedge dz_2\wedge d\bar z_1\wedge d\bar z_2 \right|
= \frac{1}{\pi^2} \vol(\CC^2/O_B(\tau_0,1)^t).
$$
Here the right-hand side is the volume under the Lebesgue measure of $\CC^2$. 

First, we claim that 
$$
\vol(\CC^2/O_B(\tau_0,1)^t)=(\mathrm{Im}(\tau_0))^2\,\vol(M_2(\RR)/\sigma(O_B)),
$$
where $M_2(\RR)\simeq\RR^4$ uses the Lebesgue measure.
In fact, let $g_1,g_2,g_3, g_4$ be a $\ZZ$-basis of $\sigma(O_B)$. 
Write $g_j=\matrixx{a_j}{b_j}{c_j}{d_j}$ for $j=1,2,3,4$. 
Write $\tau_0=x+yi$.
Then we see that $O_B(\tau_0,1)^t$ is generated by $(a_j\tau_0+b_j,c_j\tau_0+d_j)^t\in \CC^2$ for $j=1,2,3,4$. 
In terms of $\CC^2\simeq \RR^4$, this becomes $(a_jx+b_j, a_jy, c_jx+d_j, c_jy)^t$  for $j=1,2,3,4$.
They form a $4\times 4$ matrix, 
and $\vol(\CC^2/O_B(\tau_0,1)^t)$ is the absolute value of the determinant of this matrix.
By row operations, this can be changed to 
$(b_j, a_j, d_j, c_j)^t$  for $j=1,2,3,4$, and the determinant is changed by a multiple $y^2$.

Second, we claim that 
$$\vol(M_2(\RR)/\sigma(O_B))=d_B.$$ 
In fact, for any full lattice $\Lambda$ of $M_2(\RR)$, we have 
$$
\vol(M_2(\RR)/\Lambda)\vol(M_2(\RR)/\Lambda^{\#})=1
$$
for the dual lattice $\Lambda^{\#}$ of $\Lambda$ under the standard pairing $(x,y)\mapsto \tr(x y^\iota)$ on $M_2(\RR)$.
Here $y^\iota$ denotes the main involution.
This can be checked by using elementary operations in $M_2(\RR)\simeq\RR^4$ to convert $\Lambda$ to the self-dual lattice $M_2(\ZZ)\simeq\ZZ^4$.
By \cite[Lemma 15.6.17]{Voi}, the dual $\sigma(O_B)^{\#}$ of the lattice $\sigma(O_B)$ satisfies $\#(\sigma(O_B)^{\#}/\sigma(O_B))=d_B^2$.
This gives the result.

As a consequence,
$$
\| dz_1\wedge dz_2\|^2_{\Fal}(\tau_0) =  
\frac{1}{\pi^2} (\mathrm{Im}(\tau_0))^2 d_B.
$$
As $\tau_0$ varies, we have 
$$
\| dz_1\wedge dz_2\|^2_{\Fal} =  
\frac{1}{\pi^2} (\mathrm{Im}(\tau))^2 d_B.
$$

Recall that Theorem \ref{analytic} gives
$$
\psi: \underline\omega_{\CA/\CH}^{\otimes 2} \lra \omega_{\CH/\CC}^{\otimes 2},\quad
(dz_1\wedge dz_2)^{\otimes 2}\longmapsto 
\frac{d_B}{(2\pi i)^2}  (d\tau)^{\otimes 2}.
$$
Then we exactly have 
$$
\| dz_1\wedge dz_2\|^2_{\Fal}
=\| \psi(dz_1\wedge dz_2)\|^2_{\Pet}.
$$
This proves the compatibility of the metrics.

\section{Modular curves} \label{sec modular curve}

If $B=M_2(\QQ)$, then the Shimura curve $X_U$ is just the usual modular curve. In this section, we consider extensions of the Hodge bundles and the Kodaira--Spencer maps to the Deligne--Rapoport compactification of the modular curves. 
Most of the treatment is similar to the case of general Shimura curves, so our exposition here will be sketchy.

\subsection{Compactification}

Assume $B=M_2(\QQ)$ and take $O_B=M_2(\ZZ)$. 
For simplicity, we assume that
$$U=U(N)=\prod_{p\nmid N}  \GL_2(\ZZ_p) \times \prod_{p| N} (1+N M_2(\ZZ_p))^\times$$
is the principal open subgroup of $\wh O_B^\times=\GL_2(\wh\ZZ)$ for some $N\geq 1$. 
As in \S\ref{sec Shimura}, we have the complex modular curve 
$$X_U=\GL_2(\QQ) \backslash \CH^{\pm}\times \GL_2(\BA_f)/U
=\GL_2(\QQ)_+ \backslash \CH\times \GL_2(\BA_f)/U.$$
Then we have the compactification
$$\widetilde X_U=\GL_2(\QQ)_+ \backslash \CH^*\times \GL_2(\BA_f)/U,$$
where $\CH^*=\CH\sqcup \BP^1(\QQ)$ is the extended upper half plane.
The result $\wt X_U$ is a proper and smooth orbifold, which is a projective curve if  $N\geq3$.

The modular curve $X_U$ has a canonical integral model $\CX_U$ over $\ZZ[1/N]$.
In fact, $\CX_U$ is a stack over $\ZZ[1/N]$ such that for any $\ZZ[1/N]$-scheme $S$, $\CX_U(S)$ is the category of pairs 
$(E,\eta)$ as follows:
\begin{enumerate}[(1)]
\item $E$ is an abelian scheme of pure relative dimension 1 over $S$;
\item $\eta: (\ZZ/N\ZZ)^2\to E(S)[N]$ is an isomorphism of groups. 
\end{enumerate}
Then $\CX_U$ is a smooth Deligne-Mumford stack over $\ZZ[1/N]$. 

Note that the definition of $\CX_U$ here is compatible with the definition in \S\ref{sec Shimura}. 
In fact, for a triple $(A, i,\bar\eta)$ in the setting of Shimura curve, we can use the usual idempotents of $M_2(\ZZ)$ to get a canonical splitting $A\simeq E\times E$.

Following \cite[IV, \S2]{DR}, the canonical compactification $\widetilde\CX_U$ of 
$\CX_U$ is a stack over $\ZZ[1/N]$ defined as follows.
For any $\ZZ[1/N]$-scheme $S$, $\wt\CX_U(S)$ is the category of pairs 
$(\wt E,\eta)$ as follows:
\begin{enumerate}[(1)]
\item $\wt E$ is a generalized elliptic curve over $S$ whose geometric fibers are either elliptic curves or N\'eron $N$-gons;
\item $\eta: (\ZZ/N\ZZ)^2\to \wt E^\circ(S)[N]$ is an isomorphism of groups. 
\end{enumerate}
Here $\wt E^\circ$ denotes the smooth locus of $\wt E$ over $S$, which is a smooth group scheme over $S$.

It turns out that $\widetilde\CX_U$ is a proper and smooth Deligne--Mumford stack over $\ZZ[1/N]$, and it is a scheme if $N\geq3$. The natural morphism
$\CX_U\to \widetilde\CX_U$ is open, and the complement 
$$\widetilde\CX_U^\infty:=\widetilde\CX_U\setminus \CX_U$$ 
is finite and \'etale over $\ZZ[1/N]$.

\subsection{Kodaira--Spencer map}

Denote by $\pi:\CE\to \CX_U$ the universal elliptic curve, and denote by $\wt\pi:\wt\CE\to \wt\CX_U$ the universal generalized elliptic curve.
Denote by $\epsilon: \CX_U\to \CE$ and $\wt\epsilon: \wt\CX_U\to \wt\CE$ the identity sections.

The \emph{Hodge bundle} $\underline\omega_{\CE}$ over $\CX_U$ is defined as 
$$
\underline\omega_{\CE}:=\underline\Omega_{\CE}:=\epsilon^* \Omega_{\CE/\CX_U}
\simeq \pi_*\Omega_{\CE/\CX_U}.
$$
The \emph{Hodge bundle} $\underline\omega_{\wt\CE}$ over $\wt\CX_U$ is defined as 
$$
\underline\omega_{\wt\CE}:=\underline\Omega_{\wt\CE}:=\epsilon^* \Omega_{\wt\CE/\wt\CX_U}
\simeq \wt\pi_*\Omega_{\wt\CE/\wt\CX_U}.
$$

Endow the Hodge bundle $\underline\omega_{\CE}$ over $\CX_U$  with the 
\emph{Faltings metric}
$\|\cdot\|_\Fal$  as follows.
For any point $x\in \CX_U(\CC)$, and any section $\alpha\in \underline\omega_{\CE}(x)\simeq \Gamma(\CE_x,\omega_{\CE_x/\CC})$, the Faltings metric is defined by
$$
\|\alpha\|_\Fal^2:=\frac{1}{2\pi}\left|\int_{\CE_x(\CC)} \alpha\wedge\overline\alpha \right|.
$$
Here $\CE_x$ is the fiber of $\CE$ above $x$, and
$\alpha$ is viewed as a holomorphic 1-form over $\CE_x$ via the canonical isomorphism
$\underline\omega_{\CE}(x)\simeq \Gamma(\CE_x,\omega_{\CE_x/\CC})$.

Endow the relative dualizing sheaf $\omega_{\CX_U/\ZZ[1/n]}$ over $\CX_U$ with the 
\emph{Petersson metric}
$\|\cdot\|_\Pet$ as follows.
Via the complex uniformization of $\CH$ to every connected component of $\CX_U(\CC)$, the Petersson metric is defined by 
$$
\|d\tau\|_\Pet:=2\,\Im(\tau). 
$$
Here $\tau$ is the usual coordinate function of $\CH\subset\CC$.

The Kodaira--Spencer map is defined similarly. 
Start with the exact sequence
$$
0\lra \pi^*\Omega_{\CX_U/\ZZ[1/n]}\lra \Omega_{\CE/\ZZ[1/n]}\lra \Omega_{\CE/\CX_U}\lra 0.
$$
Apply derived functors of $\pi_*$.
It gives a connecting morphism 
$$
\phi_0:\pi_*\Omega_{\CE/\CX_U} \lra R^1\pi_*(\pi^*\Omega_{\CX_U/\ZZ[1/n]}).
$$
This is the \emph{Kodaira--Spencer map}. 
There are canonical isomorphisms 
$$R^1\pi_*(\pi^*\Omega_{\CX_U/\ZZ[1/n]}) \lra 
R^1\pi_*\CO_{\CE}\otimes \Omega_{\CX_U/\ZZ[1/n]}\lra 
\Lie(\CE^t)\otimes \Omega_{\CX_U/\ZZ[1/n]}
\lra \underline\Omega_{\CE^t}^\vee\otimes \Omega_{\CX_U/\ZZ[1/n]},$$ 
where $\CE^t\to\CX_U$ denotes the dual abelian scheme of $\CE\to \CX_U$.
Then the {Kodaira--Spencer map} is also written as 
$$
\phi_1:\underline\Omega_{\CE} \lra \underline\Omega_{\CE^t}^\vee\otimes \Omega_{\CX_U/\ZZ[1/n]}.
$$
Via the canonical principal polarization $\lambda:\CE\to \CE^t$ induced by line bundle associated to the identity section, the map induces a morphism 
$$
\phi_2:\underline\omega_{\CE}^{\otimes 2} \lra  \omega_{\CX_U/\ZZ[1/n]}.
$$

Note that there is a natural identification $i_{\rm can}:\CE\simeq \CE^t$ coming from identification of line bundles with divisors, but it turns out that  
$\lambda=-i_{\rm can}$. 
So switching between $\lambda$ and $i_{\rm can}$ changes $\phi_2$ by a negative sign.

As an easier version of Theorem \ref{deformation}, deformation theory simply implies that this morphism is an isomorphism over $\CX_U$. 
To extend the morphism to the compactifications, we have the following result, which is parallel to Theorem \ref{KS}.

\begin{thm}\label{KS2}
The canonical isomorphism 
$\phi_2:\underline\omega_{\CE}^{\otimes 2}\to \omega_{\CX_U/\ZZ[1/n]}$ 
induces an isomorphism
$$
\wt\phi_2:\underline\omega_{\wt\CE}^{\otimes 2}\lra
\omega_{\wt\CX_U/\ZZ[1/n]}(\wt\CX_U^\infty).
$$
Moreover, under $\phi_2$, we have $\|\cdot\|_\Fal^2=\|\cdot\|_\Pet$ over $\CX_U(\CC)$.
\end{thm}

The first statement is in \cite[VI, \S4.5]{DR}, a consequence of the calculation of deformation in \cite[III,\S1]{DR}.
The second statement is proved similarly to Theorem \ref{KS}, and it is actually slightly easier. 
We give some main ingredients below.

It suffices to compute the Kodaira--Spencer map of the universal elliptic curve
 $\pi:\CE\to \CH$
given by 
$$\CE= (\CH\times \CC)/\ZZ^2,$$ 
where 
$\ZZ^2$ acts on $\CH\times \CC$ by 
$$(a,b)\circ (\tau, z)  =(\tau, z+a \tau+b ).$$
For each $\tau\in \CH$, we have a canonical uniformization
$$\CE_\tau= \CC/\Lambda_\tau,\quad \Lambda_\tau=\tau\ZZ+\ZZ.$$
The complex torus $\CE_\tau$ has a canonical positive Riemann form 
$$
E: \Lambda_\tau \times \Lambda_\tau \lra \ZZ, \qquad 
E(a\tau+b,a'\tau+b')=ab'-a'b.
$$

Similar to Theorem \ref{analytic}, we can prove that
the Kodaira--Spencer map
$$
\phi:\pi_* \Omega_{\CE/\CH} \lra \Lie(\CE/\CH)\otimes \Omega_{\CH/\CC}
$$
gives
$$
\phi(dz)=\frac{i}{2\pi}  \frac{\partial}{\partial z}\otimes d\tau.
$$ 
Then the induced map 
$$\phi_2:\underline\omega_{\CE}^{\otimes 2}\lra \omega_{\CH/\CC}$$ 
gives 
$$
\phi((dz)^{\otimes 2})=\frac{i}{2\pi}  d\tau.
$$ 
This implies the compatibility of the metrics in Theorem \ref{KS2}.

\subsection{Coarse moduli schemes}

In \S\ref{sec coarse}, we have introduced the {coarse moduli scheme} $\CX_U^\cs$ associated to $\CX_U$. Denote by $\wt\CX_U^\cs$ the \emph{coarse moduli scheme} associated to $\wt\CX_U$.
This is also constructed by the quotient process in \S\ref{sec coarse}.
See also \cite[VI]{DR} for more properties of it.
Similar to Shimura curves, $\wt\CX_U^\cs$ is a 2-dimensional $\QQ$-factorial normal scheme, flat and projective over $\ZZ[1/N]$. 

Similar to the situation in \S\ref{sec coarse}, the relative dualizing sheaf
of the regular locus of $\wt\CX_U^\cs$ over $\ZZN$
extends to a unique $\QQ$-line bundle over $\wt\CX_U^\cs$.
Denote this extension by $\omega_{\wt\CX_U^\cs/\ZZ[1/N]}$, and call it the 
\emph{relative dualizing sheaf} of $\wt\CX_U^\cs$ over $\ZZ[1/N]$.

In \S\ref{sec coarse}, we have defined the Hodge bundle $\CL_U$ over $\CX_U^\cs$. 
The \emph{Hodge bundle} $\wt\CL_U$ over $\wt \CX_U^\cs$ is a $\QQ$-line bundle over 
$\wt\CX_U^\cs$ defined by
$$
\wt\CL_U:= \omega_{\wt\CX_{U}^\cs/\ZZ[1/N]} \otimes 
\left(\bigotimes_{P\in \wt\CX^\cs_{U,\QQ}} \CO_{\CX_U^\cs}(\CP)^{\otimes (1-e_P^{-1})}\right).
$$
Here the summation is through closed points $P$ on the generic fiber $\wt\CX^\cs_{U,\QQ}$ of $\wt\CX^\cs_{U}$ over $\QQ$, 
and $\CP$ is the Zariski closure of $P$ in $\wt\CX_U^\cs$.
If $P$ is a cusp, set $e_P=\infty$ and $1-e_P^{-1}=1$; 
if $P$ is not a cusp, then $e_P$ is the ramification index of $P$ in the map $\CX_{U(N')}\to \CX_{U}^\cs$ 
for any open subgroup $U(N')\subset U$ with $N'\geq3$.
One can check that $e_P$ does not depend on the choice of $N'$, and that 
$e_P$ is also equal to the ramification index by the uniformization map from $\CH$ to connected components of $\CX_U^\cs(\CC)$.

If $\CX_{U}$ is already a scheme,
then we simply have 
$$\wt\CX_{U}^\cs=\wt\CX_{U},\qquad\wt\CL_{U}= \omega_{\wt\CX_{U}/\ZZ[1/N]}(\wt\CX_{U}^{\infty}).$$
Here $\wt\CX_{U}^{\infty}=\wt\CX_{U}\setminus \CX_{U}$
is the locus of cusps, endowed with the reduced scheme structure.
The following result, a counterpart of Lemma \ref{compatibility}, justifies the definition involving the ramification indices.

\begin{lem} \label{compatibility2}
\begin{itemize}
\item[(1)]
Let $U'=U(N')$ be an open subgroup of $U$ for some $N'\geq3$.
Let $\wt\pi_{U',U}:\wt\CX_{U'}\to \wt\CX_{U}^\cs$ be the natural morphism.
Then there is a canonical isomorphism
$$
\wt\pi_{U',U}^*\wt\CL_{U} \lra \wt\CL_{U'}
$$
of $\QQ$-line bundles over $\wt\CX_{U'}$.
\item[(2)]
Denote by $\tilde f:\wt\CX_{U}\to \wt\CX_{U}^\cs$ the canonical morphism.
Then there is a canonical isomorphism
$$
\tilde f^*\CL_{U} \lra \omega_{\wt\CX_U/\ZZ[1/n]}(\wt\CX_U^\infty).
$$
of $\QQ$-line bundles over $\wt\CX_{U}$.
\end{itemize}
\end{lem}

\begin{proof}
The proof is similar to Lemma \ref{compatibility}, noting that the Hurwitz formula also matches multiplicities at cusps. 
\end{proof}

Recall that in \S\ref{sec coarse}, the Petersson metric \emph{Petersson metric} 
$\|\cdot\|_\Pet$ of $\CL_U$ is defined such that its pull-back to $\CH$ is given by
$$\|d\tau\|_{\mathrm{Pet}}=2\, \Im(\tau),$$
where $\tau$ is the standard coordinate function on $\CH\subset \CC$. 
Then the isomorphisms in Lemma \ref{compatibility2} are isometries.

Up to now, we have introduced the metric 
$\|\cdot\|_\Fal$ of $\underline\omega_{\CE}$ over $\CX_U$,
the metric $\|\cdot\|_\Pet$ of $\omega_{\CX_U/\ZZn}$ over $\CX_U$, and the
metric $\|\cdot\|_\Pet$ of $\CL_U$ over $\CX_U^\cs$. 
Considering the behavior of the metics at the cusps, they are not smooth, but have logarithmic singularity. 
This lies in the framework of Bost \cite{Bos} or K\"uhn \cite{Kuh}, by which their arithmetic intersection numbers are still defined.

In the case $U=\GL_2(\wh\ZZ)$, we have $\wt\CX_U^\cs\simeq \BP^1_\ZZ$ via the $j$-function. See \cite[VI, \S 1]{DR}.
Finally, we have the following counterpart of Theorem \ref{height}.

\begin{thm}\label{height2}
Assume  
$U=\GL_2(\wh\ZZ)$.
Then the normalized arithmetic intersection numbers satisfy
$$
\frac{\widehat\deg(\hat c_1(\wt\CL_U, \|\cdot\|_\Pet)^2)}{2\deg(\wt\CL_{U,\QQ})}
=\frac{\widehat\deg(\hat c_1(\omega_{\wt\CX_U/\ZZ}(\wt\CX_U^\infty), \|\cdot\|_\Pet)^2)}{2\deg(\omega_{\wt\CX_{U,\QQ}/\QQ}(\wt\CX_{U,\QQ}^\infty))}
=\frac{\widehat\deg(\hat c_1(\underline\omega_{\wt\CE}, \|\cdot\|_\Fal)^2)}{\deg(\underline\omega_{\wt\CE,\QQ})}.
$$
\end{thm}

\begin{remark}
Note that \cite[Theorem 1.1]{Yua} computes the first term of the theorem (as a special case);
\cite[Theorem 6.1]{Kuh} computes the numerator of the third term of the theorem.
The theorem asserts that these two formulas are compatible. 
\end{remark}

\

\noindent \small{Beijing International Center for Mathematical Research, Peking University, Beijing 100871, China}

\noindent \small{\it Email: yxy@bicmr.pku.edu.cn}


\begin{thebibliography}{[AB]}

\bibitem[BC]{BC}
J. -F. Boutot,  H. Carayol, Uniformisation $p$-adique des courbes de Shimura: les th\'eor\`emes de \v Cerednik et de Drinfel'd. Ast\'erisque No. 196-197 (1991), 7, 45--158 (1992). 

\bibitem[Bos]{Bos}
J. -B. Bost, 
\emph{Potential theory and Lefschetz theorems for arithmetic surfaces}, Ann. Sci. \'Ecole Norm. Sup. 32 (1999), 241--312.

\bibitem[Bou]{Bou}
J. F. Boutot, 
Le probleme de modules en inegale caracterstique in Varietes de Shimura et
fonctions L, Publ. Math. Univ. Paris VII 6, Universit\'e de Paris VII, U.E.R. de Math\'ematiques,
Paris, 1979, 43-62.

\bibitem[Buz]{Buz}
K. Buzzard, 
Integral models of certain Shimura curves.
Duke Math. J. 87 (1997), no. 3, 591--612.

\bibitem[DR]{DR}
P. Deligne, M. Rapoport, Les sch\'emas de modules de courbes elliptiques. Modular Functions of One Variable II, Proc. Internat. Summer School, Univ. Antwerp 1972, Lect. Notes Math. 349 (1973), 143--316.

\bibitem[KRY1]{KRY1}
S. Kudla, M. Rapoport, T. Yang,
{Derivatives of Eisenstein series and
Faltings heights}, Compos. Math. 140 (2004), no. 4, 887--951.

\bibitem[KRY2]{KRY2}{
S. Kudla, M. Rapoport, and T. Yang, 
{Modular forms and special cycles on Shimura curves},
Annals of Math. Studies series, vol 161, Princeton Univ. Publ., 2006.}



\bibitem[Kuh]{Kuh}
U. K\"uhn, {Generalized arithmetic intersection numbers},
J. reine angew. Math. 534 (2001), 209--236. 


\bibitem[Lan]{Lan}
K. Lan, {Arithmetic compactifications of PEL-type Shimura varieties}, London Mathematical Society Monographs Series, 36. Princeton University Press, Princeton, NJ, 2013. 


\bibitem[Mil]{Mil}
J. S. Milne,  {Points on Shimura varieties mod $p$}, Automorphic forms, representations and L-functions (Proc. Sympos. Pure Math., Oregon State Univ., Corvallis, Ore., 1977), Part 2, pp. 165--184, Proc. Sympos. Pure Math., XXXIII, Amer. Math. Soc., Providence, R.I., 1979.

\bibitem[Mum]{Mum}
D. Mumford,  
{Abelian varieties}, With appendices by C. P. Ramanujam and Yuri Manin. Corrected reprint of the second (1974) edition. Tata Institute of Fundamental Research Studies in Mathematics, 5. Published for the Tata Institute of Fundamental Research, Bombay; by Hindustan Book Agency, New Delhi, 2008.


\bibitem[Voi]{Voi}
J. Voight,  {Quaternion algebras}, Graduate Texts in Mathematics, 288. Springer, Cham, 2021.

\bibitem[Yua]{Yua}
X. Yuan, 
Modular heights of quaternionic Shimura curves, 
arXiv: 2205.13995.


\bibitem[YZ]{YZ}
X. Yuan, S. Zhang, 
{On the averaged Colmez conjecture},
Ann. of Math. (2) 187 (2018), no. 2, 533--638.

\end{thebibliography}
\end{document}